\documentclass{article}
\usepackage{amssymb}
\usepackage{amsfonts}
\usepackage{amsmath}
\usepackage{color}
\def\proof{ {\it Proof.} $\,$}
\def\sl{{\rm{SL}}}
\def\gl{{\rm{GL}}}

\def\su{{\rm{SU}}}
\def\psl{{\rm{PSL}}}
\def\psu{{\rm{PSU}}}
\def\pso{{\rm{P \Omega }}}
\def\psizo{{\rm{\Omega }}}
\def\sp{{\rm{Sp}}}
\def\psp{{\rm{PSp}}}
\def\ppsl{ ( {\rm{P}} ) {\rm{SL}}}
\def\ppsp{ ( {\rm{P}} ) {\rm{Sp}}}
\def\ppsu{ ( {\rm{P}} ) {\rm{SU}}}
\def\heart{\heartsuit}
\def\eps{\varepsilon}
\def\s{\sigma}

\def\ag{\bar{G}}
\def\at{\bar{T}}

\newtheorem{theorem}{Theorem}[section]
\newtheorem{lemma}[theorem]{Lemma}

\newtheorem{remark}[theorem]{Remark}
\newtheorem{corollary}[theorem]{Corollary}
\newtheorem{algo}[theorem]{Algorithm}

\begin{document}
\title{Construction of long root $\sl_2(q)$-subgroups in black-box groups}

\author{\c{S}\"{u}kr\"{u} Yal\c{c}\i nkaya \\
Middle East Technical University, Ankara, Turkey\\
ysukru@metu.edu.tr\\
\date{}
}
\maketitle

\begin{abstract}
We present a one sided Monte--Carlo algorithm which constructs a long root $\sl_2(q)$-subgroup in $X/O_p(X)$, where $X$ is a black-box group and $X/O_p(X)$ is a finite simple group of Lie type defined over a field of odd order $q=p^k > 3$ for some $k\geqslant 1$. Our algorithm is based on the analysis of the structure of centralizers of involutions and can be viewed as a computational version of Aschbacher's Classical Involution Theorem. We also present an algorithm which determines whether the $p$-core (or ``unipotent radical'') $O_p(X)$ of a black-box group $X$ is trivial or not, where $X/O_p(X)$ is a finite simple classical group of odd characteristic $p$. This answers a well-known question of Babai and Shalev.
\end{abstract}

\section{Introduction}

Black-box groups were introduced by Babai and Szemer\'edi \cite{babai84.229} as a generalization of permutation and matrix group algorithms in computational group theory. A \textit{black-box group} $G$ is defined to be a group equipped with a `black box' where the group elements are represented as 0-1 strings of uniform length $N$ (encoding length) and the `black box' performs the group operations which are multiplication, inversion and decision on whether a string represents a trivial element in the group. %We have an upper bound $|G|\leqslant 2^N$. For example, if $G$ is a group of Lie type of rank $n$ defined over a field of size $q$, then $O(N)=n^2\log q$.

A black-box group $G$ is specified as $G=\langle S\rangle$ for some set $S$ of elements of $G$ and to construct a black-box subgroup means to construct some generators for this subgroup. 

If $N$ is the encoding length of a black-box group $G$, then we have $|G|\leqslant 2^N$. Therefore, if $G$ is a group of Lie type of rank $n$ defined over a field of size $q$, then $|G|>q^{n^2}$, so $O(N)=n^2\log q$.

A \textit{black-box group algorithm} is defined to be  an algorithm which does not use specific properties of the representation of the given group or particulars of how the group operations are performed \cite{seress2003}. 
There are mainly two types of black-box group algorithms: Monte-Carlo and Las Vegas.  A Monte-Carlo algorithm is a randomized algorithm which may produce an incorrect output with some probability of error $\epsilon$ controlled by the user, whereas a Las Vegas algorithm never produces an incorrect answer but may report failure where the probability of the failure is less than some specified value $\epsilon$.  For a thorough discussion of randomized algorithms, see \cite{babai97.1}. The complexity of a polynomial time black-box group algorithm is of the form $O(N^c\log(1/\epsilon))$, where $N$ is the encoding length and $c$ is a constant. 

An important component of a black-box group algorithm is the construction of uniformly distributed random elements in a given group. Although not convenient for practical purposes, there is a polynomial time Monte-Carlo algorithm producing ``nearly'' uniformly distributed random elements \cite{babai91.164}. A more practical solution is the ``product replacement algorithm'' \cite{celler95.4931}, see also  \cite{pak.01.476, pak01.301}.

The main goal in this subject is to design efficient recognition algorithms of a given group. In \cite{babai99.30}, Babai and Beals proposed an approach for the construction of the composition series of black-box groups. The major component of this approach is the development of fast recognition algorithms for the finite simple groups. The black-box recognition of symmetric and alternating groups is presented in \cite{beals03.2097, bratus00.33}, and the Lie type groups are handled in the papers \cite{babai02.383, brooksbank03.162, brooksbank08.885, brooksbank01.95, brooksbank06.256, kantor01.168}. The algorithm in \cite{babai02.383} computes the standard name of a given finite simple black-box group of Lie type whereas the algorithms in \cite{kantor01.168} recognize black-box classical groups constructively, namely, if successful, the algorithm constructs an isomorphism  between the given black-box group and its standard copy. However, these constructive recognition algorithms are not polynomial in the size of the input; their complexity contain a factor of $q$ (the size of the field), whereas $N$ can contain only factors of $\log q$. Later, by assuming a procedure for the constructive recognition of $\sl_2(q)$, polynomial time constructive recognition algorithms were obtained in a series of papers \cite{brooksbank03.162, brooksbank08.885, brooksbank01.95, brooksbank06.256}.

Another important part of the construction of the composition series is to decide whether the $p$-core (or ``unipotent radical'') $O_p(X)$ of a black-box group $X$ is trivial or not, where $X/O_p(X)$ is a simple group of Lie type of characteristic $p$. It is proved in \cite{babai01.39} that if $X/O_p(X)$ is a simple unisingular group of Lie type (see Section \ref{p-core}), then random search in $X$ works effectively to decide whether $O_p(X)$ is trivial or not. Recall that a finite simple group $G$ of Lie type of characteristic $p$ is called unisingular if every nontrivial $G$-module $M$ of characteristic $p$ has the property that every element of $G$ has a nonzero fixed point in its action on $M$. The unisingular simple groups of Lie type are classified in \cite{guralnick03.271} and they constitute a small list among all simple groups of Lie type. However, if $X/O_p(X)$ is not unisingular, then the proportion of elements whose orders are multiples of $p$ is $O(1/q)$, where $q$ is the size of the underlying field \cite{babai01.39, guralnick01.169}. Moreover, deciding whether $O_p(X)$ is trivial or not is a harder problem \cite{babai01.39}. 

In \cite{suko01}, the author proposed and briefly outlined a uniform approach recognizing black-box groups $X$ where $X/O_p(X)$ is a Lie type group of odd characteristic $p$. The main idea is to construct all root $\sl_2(q)$-subgroups in a Lie type group corresponding to the nodes in the extended Dynkin diagram of the corresponding algebraic group. The approach is based on recursive construction of centralizers of involutions in black box groups \cite{altseimer01.1,  borovik02.7, bray00.241}. The reader is referred to \cite{suko01} for a detailed discussion of this approach together with the previous algorithms recognizing the black-box groups of Lie type \cite{babai02.383, brooksbank03.162,brooksbank08.885, brooksbank01.95, brooksbank06.256, kantor01.168}. 

The object of this paper is to complete the first part of this project, that is, to present and justify an algorithm which constructs a subgroup corresponding to a long root $\sl_2(q)$-subgroup in a black-box group $X$ where $X/O_p(X)$ is a Lie type group of odd characteristic  $p$ and decides whether the $p$-core (or ``unipotent radical'') $O_p(X)$ is trivial or not.

The main result is the following.

\begin{theorem}\label{main1}
Let $X$ be a black-box group and $\bar{X}=X/O_p(X)$ isomorphic to a finite simple group of Lie type defined over a field of odd order $q = p^k > 3$ for some $k\geqslant 1$.  Assume that $\bar{X} \ncong \psl_2(q)$ and $\bar{X} \ncong \/ ^2G_2(q)$, then there is a polynomial time Monte-Carlo algorithm which constructs a subgroup $K$ such that $K/O_p(K)$ is a long root $\sl_2(q)$-subgroup in $\bar{X}$.
\end{theorem}

If $X/O_p(X)$ is a Lie type group defined over a field of odd characteristic $p$ with non-trivial center, then $Z(X)$ can be constructed by a polynomial time Monte--Carlo algorithm \cite{babai01.39}. Therefore the algorithm in Theorem \ref{main1} can be extended to all quasi-simple groups of Lie type over a field of odd order $q >3$ except for $X/O_p(X) \ncong \psl_2(q)$ or $^2G_2(q)$.

In the case of a black-box group $X$ where $\bar{X} = X/O_p(X)$ is isomorphic to $\psl_2(q)$ or $^2G_2(q)$, there is no subgroup in $\bar{X}$ isomorphic to $\sl_2(q)$. Therefore we exclude these groups in Theorem \ref{main1}. 

In a later publication \cite{suko03}, we extend this algorithm to construct all subgroups $K\leqslant X$ where $K/O_p(K)$ correspond to the root $\sl_2(q)$-subgroups appearing as a node in the extended Dynkin diagram of $X/O_p(X)$. In other words, we construct the Curtis-Tits and Phan systems \cite{curtis65.174, phan77.67, timmesfeld04.38} for the black-box groups of Lie type of odd characteristic $p$. As discussed in \cite[Section 2]{suko01}, we note here that this approach can be seen as a computational version of Aschbacher's ``Classical Involution Theorem" \cite{aschbacher77.353} which plays a prominent role in the classification of the finite simple groups. We also note that this procedure is not a constructive recognition algorithm.

The next result allows us to decide whether a given subgroup is a long root $\sl_2(q)$-subgroup in a finite simple group of Lie type defined over a field of odd order $q>3$. 

\begin{theorem}\label{main2}
Let $K$ be a black-box subgroup in a finite simple black-box group of Lie type defined over a field of odd order $q>3$ isomorphic to $\ppsl_2(q^k)$ for some $k\geqslant 1$. Then there is a polynomial time Monte-Carlo algorithm which decides whether $K$ is a long root $\sl_2(q)$-subgroup.
\end{theorem}

One of the most important applications of Theorem \ref{main1} is a polynomial time Monte-Carlo algorithm which determines whether $O_p(X) \neq 1$ answering a well-known question of Babai and Shalev \cite{babai01.39}. We prove the following theorem.

\begin{theorem}\label{maincor1}
Let ${X}$ be a black-box group with the property that $X/O_p(X)$ is a simple classical group of odd characteristic $p$.  Then there is a polynomial time Monte-Carlo algorithm which decides whether $O_p(X) \ne 1$, and, if $O_p(X) \ne 1$, we can find a non-trivial element of $O_p(X)$.
\end{theorem}

Combining Theorem \ref{maincor1} with the results in \cite{altseimer01.1, babai99.30, babai02.383, kantor02.370} we obtain:

\begin{theorem}
Finite simple classical groups of known odd characteristic can be recognized in Monte--Carlo polynomial time among all black-box groups.  
\end{theorem}

The black-box group operations allows us to work only with the multiplication table of the given group. Therefore it is almost impossible, for example, in big matrix groups, to get information about the group without an additional oracle. One can solve many problems by introducing an \textit{order oracle} with which we can find the orders of elements. In this paper, we do not attempt to find the exact orders of elements. Instead we work with a milder assumption that we are given a computationally feasible global exponent $E$ for $X$ as an input, that is, a reasonably sized natural number $E$ such that $x^E = 1$ for all $x \in X$. Notice that one can set $E=|\gl_n(q)|$ for $n \times n$ matrix groups. 

The characteristic of the underlying field for the black-box groups of Lie type can be computed by using an order oracle \cite{kantor02.370}. A more recent algorithm computing the characteristic, which makes use of an order oracle and the construction of centralizers of involutions, is presented in \cite{liebeck07.741}. In our algorithms we assume that the characteristic $p$ of the underlying field is given as an input.

In this paper we are mainly interested in the case where the order of the base field is large. Note that the algorithms also work for groups over small fields provided that solvable $2$-components do not occur in the centralizers of involutions, or equivalently, that $q>3$. However there are better algorithms where $X/O_p(X)$ can be recognized constructively when the field of definition is small \cite{kantor01.168}.

The algorithms presented in this paper have been tested extensively in GAP on various groups over a variety of fields. The implementation is discussed in Section \ref{gap}.

\section{Background}
In this section, we recall some basic properties of groups of Lie type that we use frequently in the present paper. We refer the reader to  \cite{carter1972, carter1985, gorenstein1998} for a complete description of groups of Lie type.

%of maximal tori and centralizers of involutions in finite groups of Lie type.  For a brief summary of the properties of root subgroups and root $\sl_2(q)$-subgroups that we need in this paper we refer the reader to \cite{suko01}. Standard references for a complete description of groups of Lie type are \cite{carter1972, carter1985, gorenstein1998}.

Throughout this section we set the following notation. Let  $\ag$ denote a connected simple algebraic group over an algebraically closed field of characteristic $p$, $\at$ a maximal torus of $\ag$, $\bar{B}$ a Borel subgroup containing $\at$ and $\bar{\Sigma}$ be the corresponding root system. Let $\bar{N} = N_{\ag}(\at)$ and $W=\bar{N}/\at$ the Weyl group of $\ag$. Let $\s$ be a Frobenius endomorphism of $\ag$ and $\ag_\s$ the fixed point subgroup of $\ag$ under $\s$.  The subgroup $\ag_\s$ is finite, and we denote $G=O^{p^\prime}(\ag_\s)$. 

%Throughout this section, we set a maximal torus $\at$ and Borel subgroup $\bar{B}$ of $\ag$ which are $\s$-invariant and $\at \leqslant \bar{B}$ \cite[2.1.6]{gorenstein1998}. 

\subsection{Maximal tori}\label{maximal_tori}
It is known that there exists a maximal torus $\at$ and Borel subgroup $\bar{B}$ of $\ag$ which are $\s$-invariant and $\at \leqslant \bar{B}$, see for example \cite[2.1.6]{gorenstein1998}. The subgroup of the form $G\cap \at$ for some $\s$-invariant maximal torus $\at$ of $\ag$ is called a \textit{maximal torus} of $G$.

The fundamental result about maximal tori in finite groups of Lie type is that the set of $G$-orbits on the set of $\s$-invariant maximal tori of $\ag$ is in bijective correspondence with $H^1(\s,W)$ \cite[2.1.7]{gorenstein1998}, where $H^1(\s,G)$ is the set of equivalence classes of $G$ (or $\s$-conjugacy class of $G$) under the relation $\sim$ defined by
\begin{equation}\label{max-torus}
x \sim y \,\, {\rm{if \, and\, only \, if}} \,\, y=gxg^{-\s} \, {\rm{for \, some}}\,\, g \in G.
\end{equation}
In particular, if $\s$ fixes every element of $W$, which is the case if $G$ is untwisted, then $H^1(\s,W)$ corresponds to the set of conjugacy classes of $W$. If a torus $T\leqslant G$ corresponds to a $\s$-conjugacy class of $W$ containing $w$, then $T$ is called a torus \textit{twisted by $w$}. 

The set $C_W(w,\s) = \{ x \in W \mid   w=xwx^{-\s} \}$ is a subgroup of $W$. By \cite[3.3.6]{carter1985}, we have
\begin{equation}\label{toriat-q}
|\bar{N}_\s /\at_\s | \cong |C_W(w,\s)|,
\end{equation}
where $\bar{N} = N_{\ag}(\at)$. Hence, we obtain an important result that the sizes of the conjugacy classes of tori depend only on the Weyl group. Notice that if $G$ is an untwisted group of Lie type, then $C_W(w,\s)=C_W(w)$.

An element which belongs to a torus is called a \textit{semisimple element}, and $t \in G$ is called \textit{regular semisimple} if $\bar{T} = C_{\ag}(t)$ is a $\s$-invariant maximal torus of $\ag$.

\subsection{Root $\sl_2(q)$-subgroups}\label{root-sl2}

For each root $r \in \bar{\Sigma}$, there exists a $\bar{T}$-root subgroup of $\bar{G}$. If $\at$ is $\sigma$-invariant, then the map $\s$ permutes these root subgroups and induces an isometry on the Euclidean space $\mathbb{R}\bar{\Sigma}$ spanned by $\bar{\Sigma}$. Let $\Delta$ be a $\langle \s \rangle$-orbit of a root subgroup of $\bar{G}$, then the subgroup $O^{p^{\prime}}(\langle \Delta \rangle_\s)$ is called a $T$-root subgroup of $G$, where $T=\bar{T}_\s$ \cite{seitz83.153}. The root system of the finite group $G$ is obtained by taking the fixed points of the isometry induced from $\s$ on $\mathbb{R}\bar{\Sigma}$ \cite[Section 2.3]{gorenstein1998}, and $G$ is generated by the corresponding root subgroups. A root subgroup is called a long or short root subgroup, if the corresponding root is long or short, respectively. We refer the reader to \cite[Table 2.4]{gorenstein1998} for a complete description of the structure of the root subgroups in finite groups of Lie type.

%\begin{table}[h]
%\caption{The structure of root subgroups in $G$ \cite[Table 2.4]{gorenstein1998}. Here $E_{q^i}$ is an elementary abelian group of order $q^i$.}\label{roottwisted}
%\begin{tabular}{ccccccc}
%Type                                    && $r$            &&&  Remarks \\[.8ex] \hline 
%Untwisted                               &&    both        &&&  $X_r \cong E_q$  \\[.8ex]
%Twisted except $^2A_{2n}(q)$            &&    long        &&&  $X_r \cong E_q$  \\[.8ex]
%Twisted except $^3D_4(q)$               &&    short       &&&  $X_r \cong E_{q^2}$      \\[.8ex]
%$^3D_4(q)$                              &&    short       &&&  $X_r \cong E_{q^3}$       \\[.8ex]
%$^2A_{2n}(q)$                           &&    long        &&&  $|X_r| = q^3$ and $Z(X_r) \cong E_q$ \\[.8ex]
%$^2G_2(q)$                              &&                &&&  $|X_r| = q^6$ and $Z(X_r) \cong E_{q^2}$\\[.8ex]
%\end{tabular}
%\end{table}

Let $\Sigma=\{ r_1, \ldots , r_n\}$ and $X_{r_1}, \ldots, X_{r_n}$ be the corresponding root subgroups. Set $M_i = \langle X_{r_i}, X_{-r_i} \rangle$, $Z_i = Z(X_{r_i})$ and $K_i = \langle Z_i, Z_{-i}\rangle$. Then $X_{r_i}$ is a Sylow $p$-subgroup of $M_i$ and $Z_i$ is a Sylow $p$-subgroup of $K_i$. The subgroup $K_i\leqslant G$ is called \textit{long or short root $\sl_2(q)$-subgroup}, if the corresponding root $r_i \in \Sigma$ is a long or short, respectively. Here, $q$ is the order of the center of a long root subgroup of $G$.
In this case we say that $G$ is defined over a field of order $q$.

\begin{theorem}{\rm (\cite[Theorem 14.5]{aschbacher77.353})}\label{long-root-sl2}
Let $G$ be a finite simple group of Lie type defined over a field of odd order $q>3$ different from $\psl_2(q)$ and $^2G_2(q)$. With the above notation, let $r_i$ be a long root, $K=K_i$, and $\langle z \rangle = Z(K)$. Then
\begin{itemize}
\item[{\rm (1)}] $K\cong \sl_2(q)$.
\item[{\rm (2)}] $O^{p^\prime}(N_G(K)) = KL$, where $[K,L]=1$ and $L$ is the Levi factor of the parabolic subgroup $N_G(Z_i)$.
\item[{\rm (3)}] $K \unlhd C_G(z)$.
\end{itemize} 
\end{theorem}

\subsection{The structure of the centralizers of involutions}\label{centralizers}
In this section we summarize the structure of the centralizers of involutions in finite simple groups of Lie type of odd characteristic $p$ in Table \ref{cent-table}, which is extracted from \cite{gorenstein1998} for the convenience of the reader. We refer the reader to \cite[\S4.5, \S4.6]{gorenstein1998} for a complete description.

We use the following notation in Table \ref{cent-table}. For $\eps= \pm$, we write $\ppsl_n^\eps(q)$ to denote $\ppsl_n(q)$ and $\ppsu_n(q)$ when $\eps=+$ and $\eps = -$, respectively. Similarly $E_6^\eps(q)$ denotes the simple group $E_6(q)$, if $\eps =+$, and $^2E_6(q)$, if $\eps =-$. Moreover we denote the universal version of $E_6^\eps(q)$ by $\hat{E}_6^\eps(q)$ which is a central extension of the simple group $E_6^\eps(q)$ by an abelian group of order $(3, q-\eps)$.

We write $G_1 \circ_n G_2$ which is meant to be $(G_1 \times G_2)/N$ for some cyclic group $N$ of order $n$ intersecting with $G_1$ and $G_2$ trivially. 

Let $G$ be a finite group and $Z(G)$ is cyclic. Then we write $\frac{1}{m}G$ to denote the quotient group $G/Y$, where $Y \leqslant Z(G)$ and $|Y|=m$.
%classification book, volume 3, page 148.
Note that the center of $G\cong {\rm Spin}_{2n}^+(q)$ is an elementary abelian 2-group of order $4$ when $n$ is even. Therefore $\frac{1}{2}{\rm Spin}_{2n}^+(q)$ is not uniquely defined for $n$ even and we define it as follows. There is an involution $z \in Z(G)$ such that $G/\langle z \rangle \cong {\rm SO}_{2n}^+(q)$. For the other involutions $z_1,\, z_2 \in Z(G) \backslash \{z\}$, we have $G/\langle z_1 \rangle \cong G/\langle z_2 \rangle$ which is not isomorphic to ${\rm SO}_{2n}^+(q)$ and we denote these quotient groups as $\frac{1}{2}{\rm Spin}_{2n}^+(q)$. Notice that $\frac{1}{2}{\rm Spin}_{12}^+(q)$ and $\frac{1}{2}{\rm Spin}_{16}^+(q)$ appear as components in the centralizers of certain involutions in $E_7(q)$ and $E_8(q)$, respectively.        

The involutions in long root $\sl_2(q)$-subgroups are called \textit{classical involutions}. 
Below, we give the list of classical involutions and the semisimple socles of their centralizers.

\begin{table}[h]
\caption{Classical involutions and the semisimple socles of their centralizers.} \label{long-root}
\begin{tabular}{cccccc}
$G$                         &&  $i$  &&&  $O^{p'}(C_G(i))$ \\[.8ex] \hline
$\psl_n^\eps(q)$            && $t_2$ &&&  $\sl_2(q)\circ_2 \sl_{n-2}^\eps(q)$ \\[.8ex]
$\psp_{2n}(q)$              && $t_1$ &&&  $\sl_2(q)\circ_2 \sp_{2n-2}(q)$ \\[.8ex]
$\psizo_{2n+1}(q)$          && $t_2$ &&&  $(\sl_2(q)\circ_2 \sl_2(q))\circ_2 \psizo_{2n-3}(q)$ \\[.8ex]
$\pso_{2n}^\eps(q)$         && $t_2$ &&&  $(\sl_2(q)\circ_2 \sl_2(q))\circ_2\psizo^\varepsilon_{2n-4}(q)$ \\[.8ex]
$G_2(q)$                    && $t_1$ &&&  $\sl_2(q)\circ_2 \sl_2(q)$ \\[.8ex]
$^3D_4(q)$                  && $t_2$ &&&  $\sl_2(q)\circ_2 \sl_2(q^3)$ \\[.8ex]
$F_4(q)$                    && $t_1$ &&&  $\sl_2(q) \circ_2\sp_6(q)$ \\[.8ex]
$E_6^\eps(q)$               && $t_2$ &&&  $\sl_2(q)\circ_2 \frac{1}{(q-\eps,3)}\sl_6^\eps(q)$ \\[.8ex]
$E_7(q)$                    && $t_1$ &&&  $\sl_2(q)\circ_2 \frac{1}{2}{\rm Spin}_{12}(q)$ \\[.8ex]
$E_8(q)$                    && $t_8$ &&&  $\sl_2(q)\circ_2 E_7(q)$ \\
\end{tabular}
\end{table}

\begin{table}[p]
\caption{Centralizers of involutions in finite simple groups of Lie type of odd characteristic.} \label{cent-table}
\begin{tabular}{|c|c|c|c|c|} \hline 
$G$         &   conditions          & type    &   $O^{p^\prime}(C_G(i))$               & $|C_{C_G(i)}(L)|$  \\ [.7ex] \hline  
            &                                & $t_1$   &   $\sl_{n-1}^\eps(q)$                               & $q-\eps$         \\ [.7ex]
$\psl_n^\eps(q)$ &  $2 \leqslant k \leqslant n/2$ & $t_k$   &          $ \sl_k^\eps(q) \circ \sl_{n-k}^\eps(q)$        & $q-\eps$         \\[.7ex]

            &  $n$ even                      & $t_{n/2}^\prime$   & $\frac{1}{(n/2,q-\eps)}\sl_{n/2}(q^2)$ & $q+\eps$      \\[.7ex] \hline
                   &                     & $t_1$        & $\psizo_{2n-1}(q)$                            &$2(q-1)$ \\[.7ex]
                   &                     & $t_1^\prime$ & $\psizo_{2n-1}(q)$                            &$2(q+1)$ \\[.7ex]
                  
$\psizo_{2n+1}(q)$ &  $2\leqslant k <n $ & $t_k$            & $\psizo^+_{2k}(q) \times \psizo_{2(n-k)+1}(q)$ &  $2$ \\[.7ex]
$n\geqslant 2$     &  $2\leqslant k <n $ & $t_k^\prime$     & $\psizo^-_{2k}(q) \times \psizo_{2(n-k)+1}(q)$ &  $2$ \\[.7ex]
                   &                     & $t_n$            & $\psizo^+_{2n}(q)$                             & $2$  \\ [.7ex]
                   &                     & $t^\prime_n$     & $\psizo^-_{2n}(q)$                             & $2$  \\ [.7ex] \hline 
                   
$\psp_{2n}(q)$     & $1\leqslant k \leqslant n/2$ & $t_k$        & $\sp_{2k}(q) \circ_2 \sp_{2(n-k)}(q)$ & $2$ \\[.7ex]
$n\geqslant2$      &                              & $t_n$        & $\frac{1}{(2,n)}\sl_n(q)$           & $q-1$ \\[.7ex]
                   &                              & $t_n^\prime$ & $\frac{1}{(2,n)}\su_n(q)$           & $q+1$ \\[.7ex]  \hline
                   &                       & $t_1$        & $\psizo^\eps_{2n-2}(q)$                & $q-1$ \\[.7ex]                           &                       & $t_1^\prime$ & $\psizo^{-\eps}_{2n-2}(q)$                & $q+1$ \\[.7ex]    
${\rm P}\psizo^\eps_{2n}(q)$  & $2\leqslant k < n/2$  & $t_k$            & $\psizo^+_{2k}(q) \circ_2 \psizo^\eps_{2(n-k)}(q)$  & $2$ \\[.7ex]              
$ n\geqslant 4$            & $2\leqslant k < n/2$  & $t_k^\prime$     & $ \psizo^-_{2k}(q) \circ_2 \psizo^{-\eps}_{2(n-k)}(q)$ & $2$ \\[.7ex]   
                           & $\pso^+_{4m}(q)$  & $t_{n/2}$        & $ \psizo^+_{2m}(q) \circ_2 \psizo^+_{2m}(q)$           & $2$ \\[.7ex]                     
                           & $\pso^+_{4m}(q)$  & $t_{n/2}^\prime$ & $ \psizo^-_{2m}(q) \circ_2 \psizo^-_{2m}(q) $          & $2$ \\[.7ex]  
                           & $\pso^+_{4m}(q)$  & $t_{n-1}, \, t_n$        & $\frac{1}{2}\sl_{2m}(q)$                          & $q-1$ \\[.7ex]            
                           & $\pso^+_{4m}(q)$  & $t_{n-1}^\prime, \, t_n^\prime$ & $\frac{1}{2}\su_{2m}(q)$                           & $q+1$ \\[.7ex]
                           & $\pso^-_{4m}(q)$  & $t_{n/2}$ & $ \psizo^-_{2m}(q) \times \psizo^+_{2m}(q) $                            & $2$ \\[.7ex]
                           & $\pso^\eps_{2(2m+1)}(q)$ & $t_n$        & $\sl^\eps_{2m+1}(q)$                                     & $q-\eps$ \\[.7ex]  \hline
$^3D_4(q)$             &                        & $t_2$        & $\sl_2(q) \circ_2 \sl_2(q^3)$              & $2$\\[.7ex] \hline
$G_2(q)$               &                        & $t_1$        & $\sl_2(q)\circ_2 \sl_2(q)$                 & $2$\\[.7ex]  \hline
%$^2G_2(q)$             &                        & $t_1$        & $\psl_2(q^2)$                            & $2$\\[.7ex] \hline 
$F_4(q)$               &                        & $t_1$        & $\sl_2(q) \circ_2 \sp_6(q)$                & $2$\\[.7ex]
                       &                        & $t_4$        & ${\rm Spin}_9(q)$                              & $2$\\[.7ex] \hline 
$E_6^\varepsilon(q)$   &                        & $t_1$        & ${\rm Spin}^\varepsilon_{10}(q)$               & $q- \varepsilon $ \\[.7ex]
                       &                        & $t_2$        & $\sl_2(q) \circ_2 \frac{1}{(q-\eps,3)}\sl^\varepsilon_6(q)$    & $2$\\[.7ex] \hline 
$E_7(q)$               &                        & $t_1$        & $\sl_2(q) \circ_2 \frac{1}{2}{\rm Spin}_{12}^+(q)$ &$2$\\[.7ex]
                       &                        & $t_4$, $t_4^\prime$& $\frac{1}{(4,q-\eps)} \sl_8^\eps(q)$  & $2$\\[.7ex] 
                       &                        & $t_7$, $t_7^\prime$   &   $\hat{E}_6^\eps(q)$   & $q-\eps$\\[.7ex]   \hline         
$E_8(q)$               &                        & $t_1$        & $\frac{1}{2}{\rm Spin}_{16}^+(q)$                & $2$\\[.7ex]
                       &                        & $t_8$        & $\sl_2(q) \circ_2 E_7(q)$                  & $2$\\[.7ex] \hline
\end{tabular}
\end{table}

The following is a direct consequence of \cite[Theorem 4.2.2]{gorenstein1998} which we need in the sequel.

\begin{lemma}\label{semisocle}
Let $G$ be a finite simple group of Lie type defined over a field of odd order $q>3$ and $i \in G$ be an involution. Then the second derived subgroup $C_G(i)^{\prime \prime}$ is the semisimple socle $O^{p'}(C_G(i))$ of $C_G(i)$ which is a central product of (quasi)simple groups of Lie type in same characteristic.
\end{lemma}

Note that if $G$ is a universal version of a Lie type group defined over a field of odd order $q>3$ then $C_G(i)^\prime$ is the semisimple socle of $C_G(i)$.

Passing to the groups with a non-trivial $p$-core we have the following.

\begin{corollary}\label{semisocle2}
Let $X$ be a finite group. Assume that $X/O_p(X)$ is a finite simple group of Lie type over a field of odd order $q>3$ and $i \in X$ be an involution. Then $(C_X(i)/O_p(C_X(i)))^{\prime\prime}$ is the semisimple socle of $C_X(i)/O_p(C_X(i))$.
\end{corollary}

\section{Pairs of long root $\sl_2(q)$-subgroups}\label{pairs}
In this section, we determine the structure of the subgroups generated by randomly chosen two conjugate long root $\sl_2(q)$-subgroups in a finite simple classical group of odd characteristic. 

\begin{lemma}\label{psl}
Let $G\cong \psl_n(q)$, $n\geqslant 5$. Let $K\leqslant G$ be a long root $\sl_2(q)$-subgroup and $g \in G$ be a random element. Then $\langle K, K^g \rangle \cong \sl_4(q)$ with probability at least $1-1/q^{n-3}$.
\end{lemma} 
\proof  Let $V$ be a natural module for $\sl_n(q)$ and $V= U \oplus W$, where $K$ induces $\sl_2(q)$ on $U$ and fixes $W$. Assume that $U=\langle u_1,u_2 \rangle$. Let $g \in \sl_n(q)$ and 
$$g u_1  = a_1 u_1 + a_2 u_2 + w_1$$
$$g u_2  = b_1 u_1 + b_2 u_2 + w_2,$$
where $a_i, b_i$, $i=1,2$, are elements in the base field and $w_1, \, w_2 \in W$. Observe that the vectors $w_1$ and $w_2$ are linearly dependent with probability $1/q^{n-3}$. Therefore dim$\langle U, gU \rangle = 4$ with probability at least $1-1/q^{n-3}$.  \hfill $\Box$

The following two lemmas correspond to Lemmas 4.11 and 5.8 in \cite{kantor01.168}.

\begin{lemma}\label{psp}
Let $G\cong  \psp_{2n}(q)$, $n\geqslant 3$. Let $K \leqslant G$ be a long root $\sl_2(q)$-subgroup and $g \in G$ be a random element. Then $\langle K, K^g \rangle \cong \sp_4(q)$ with probability at least $1-1/q^{2n-2}$. 
\end{lemma}
\proof Let $V$ be the natural module for $\sp_{2n}(q)$ and $V = V_1 \bot \ldots \bot V_n$ be an orthogonal decomposition of $V$, where $V_k$ is a hyperbolic plane for each $k=1, \ldots ,n$. Assume that $V_k = \langle e_k, f_k \rangle$, where $\{e_k, f_k\}$ are hyperbolic pairs. We may assume that $K = \sp(V_1)$. Let $g \in \sp_{2n}(q)$ and $K^g= \sp(V_1^\prime)$, where $V_1^{\prime} \leqslant V$. It is clear that $V_1^\prime$ is a hyperbolic plane in $V$ and the subspace $\langle V_1, V_1^\prime \rangle$ is non-degenerate 4-space with probability 
$$ 1-A/B,$$
where $A$ is the number of hyperbolic planes intersecting with $V_1$ non-trivially and $B$ is the total number of hyperbolic planes in $V$. By the computation in \cite[Chapter 3]{artin1957}, there are $q^{2n-1}(q^{2n}-1)$ hyperbolic planes in a $2n$-dimensional symplectic space, and there are $q^{2n-1}$ hyperbolic planes  containing a fixed vector. Hence, the the probability that $\langle V_1, V_1^\prime \rangle$ is a 4-dimensional symplectic space with probability at least 
$$1-\frac{q^{2n-1}(q^2-1)}{q^{2n-1}(q^{2n}-1)} > 1-1/q^{2n-2}.$$ \hfill $\Box$

\begin{lemma}\label{psizo}
Let $G\cong  \pso^\pm(V)$ be a simple orthogonal group with ${\rm{dim}}V \geqslant 9$. Let $K \leqslant G$ be a long root $\sl_2(q)$-subgroup and $g \in G$ be a random element. Then $\langle K, K^g \rangle \cong \psizo_8^+(q)$ with probability at least $(1-1/q)^2$. 
\end{lemma}
\proof Let $V$ be the natural module for $\psizo^\pm(V)$. Then 
$$[K,V] = \langle kv-v \mid v \in V, \, k \in K \rangle$$ 
is an orthogonal $4$-space of Witt index $2$ and let $U=[K,V] = \langle e_1,e_2,f_1,f_2 \rangle$, where $\{e_i,f_i\}$ are hyperbolic pairs. 
Let $U^\prime = \langle e_1^\prime, e_2^\prime, f_1^\prime, f_2^\prime \rangle \leqslant V$ be the orthogonal $4$-space of Witt index $2$ on which $K^g$ induces $\sl_2(q)$ for $g \in \psizo^\pm(V)$. 
Again by  computation in \cite[Chapter 3]{artin1957}, there are 
\begin{center}
\begin{tabular}{rl}
 $q^{n-2}(q^{n-1}-1)$ & if $n$ is odd, \\
 $q^{n-2}(q^{n/2}-1)(q^{{n/2}-1}+1)$ & if $n$ is even and $G\cong  \pso^+(V)$, \\
 $q^{n-2}(q^{n/2}+1)(q^{{n/2}-1}-1)$ & if $n$ is even and $G\cong  \pso^-(V)$ \\ 
\end{tabular}
\end{center}
total number of hyperbolic pairs and $q^{n-2}$ hyperbolic pairs containing a fixed singular vector. By the similar computations as in Lemma \ref{psp}, the subspace $U_1 = \langle U, e_1^\prime , f_1^\prime \rangle $ is a non-degenerate 6-space with Witt index 3 and the subspace $\langle U_1 ,e_2^\prime ,f_2^\prime \rangle$ is a non-degenerate 8-space with Witt index 4 with probability at least $1-1/q$. Hence, $\langle U, gU \rangle$ is an orthogonal 8-space with Witt index 4 with probability at least $(1-1/q)^2$. \hfill $\Box$

\begin{lemma}\label{psu}
Let $G\cong \psu_n(q)$, $n\geqslant 5$. Let $K\leqslant G$ be a long root $\sl_2(q)$-subgroup and $g \in G$ be a random element. Then $\langle K, K^g \rangle \cong \su_4(q)$ with probability at least $1-1/q$.
\end{lemma}
\proof The proof is same as the proof of the previous lemmas. The computation follows from the fact that there are $q^{2n-3}(q^{n-1}-(-1)^{n-1})(q^n-(-1)^n)$ hyperbolic pairs and $q^{2n-3}$ of them contains a fixed isotropic vector in an $n$-dimensional unitary space \cite{taylor1992}. 

\begin{lemma}\label{others}
Let $G\cong  \ppsl_n(q), \, \ppsu_n(q)$, $n=2,3,4$, or $\ppsp_4(q),$ $\psizo_7(q)$, $\psizo_8^\pm(q)$, and $K$ be a long root $\sl_2(q)$-subgroup in $G$, then $\langle K, K^g \rangle = G$ with probability at least $1-1/q$.
\end{lemma}
\proof Similar to the arguments above.

\section{Construction of $C_G(i)$ in a black-box group}\label{const-bb}
In this section, we summarize the construction of the centralizers of involutions in black-box groups following \cite{borovik02.7}, see also \cite{bray00.241}.

Let $X$ be a black-box finite group having an exponent $E=2^km$ with $m$ odd. To produce an involution in $X$, we need an element $x$ of even order. Then the last non-identity element in the sequence
$$1 \neq x^{m}, \, x^{m2}, \, x^{m2^2}, \, \ldots , x^{m2^{k-1}}, x^{m2^k}=1$$ 
is an involution and denoted by ${\rm i}(x)$.

We call an element $j \in X$ of order $4$ a \textit{pseudo-involution}, if it is an involution in $X/Z(X)$ but not in $X$. Let $Y\leqslant X$, then we call an element $j \in Y$ a \textit{pseudo-involution} in $Y$, if $1 \neq j^2 \in Z(Y)$. We can produce pseudo-involutions in a black-box group $X$ in a similar manner, that is, we first produce an element of order $4$ and check whether $j^2$ commutes with the generators of $X$.

Let $i$ be an involution in $X$. Then, by \cite[Section 6]{borovik02.7}, there is a partial map $\zeta^i = \zeta^i_0 \sqcup \zeta^i_1$ defined by
\begin{eqnarray*}
\zeta^i: X & \longrightarrow &  C_X(i)\\
x & \mapsto & \left\{ \begin{array}{ll}
\zeta^i_1(x) = (ii^x)^{(m+1)/2}\cdot x^{-1} & \hbox{ if } o(ii^x) \hbox{ is odd}\\
\zeta^i_0(x) = {\rm i}(ii^x)  &  \hbox{ if } o(ii^x) \hbox{ is even.}
\end{array}\right.
\end{eqnarray*}

Here $o(x)$ is the order of the element $x\in X$. Notice that, with a given exponent $E=2^km$, we can construct $\zeta_0^i(x)$ and $\zeta_1^i(x)$ without knowing the exact order of $ii^x$.

The following theorem is the main tool in the construction of centralizers of involutions in black-box groups.

\begin{theorem} {\rm (\cite{borovik02.7})} \label{dist}
Let $X$ be a finite group and $i \in X$ be an involution. If the elements $x \in X$ are uniformly distributed and
independent in $X$, then 
\begin{enumerate}
\item the elements $\zeta^i_1(x)$ are uniformly distributed and independent in $C_X(i)$ and
\item the elements $\zeta^i_0(x)$ form a normal subset of involutions in $C_X(i)$.
\end{enumerate}
\end{theorem}

We will use both of the functions $\zeta^i_0$ and $\zeta^i_1$ to generate $C_X(i)$. It follows directly from Theorem \ref{dist} that the image of the function $\zeta^i_1$ is $C_X(i)$ and the image of $\zeta^i_0$ generates a normal subgroup in $C_X(i)$. Although the map $\zeta_1$ is a better black-box for the construction of centralizers of involutions, it turns out that the function $\zeta_0$ is sufficient for our purposes, see Section \ref{section-heart}.

\section{The heart of the centralizer}\label{section-heart}
In this section, we describe the subgroup generated by the image of the function $\zeta^i_0$ for any involution $i \in G$, where $G$ is a finite simple group of Lie type of odd characteristic. 

Let $i \in G$ be an involution. Define
$$\heart_i(G) = \langle \zeta^i_0(g) \mid g \in G \rangle.$$
Here, we use the convention that $\zeta^i_0(g)=1$, if $ii^g$ has odd order. We also assume in the sequel that $\zeta_1^i(g)=1$ when $ii^g$ has even order.

\begin{lemma}\label{zeta0-lemma}
Let $G$ be a finite group and $i \in G$ be an involution. Then the image of $\zeta^i_0$ does not contain involutions from the coset\/ $iZ(G)$.
\end{lemma}
\proof Assume that  $\zeta^i_0(g) \neq 1$ for some $g \in G$, and consider the dihedral group $D\/=\/\langle i, i^g \rangle$. Recall that $\zeta^i_0(g) = {\rm i}(ii^g) \in Z(D)$. Therefore, if $\zeta^i_0(g) = iz$, where $z \in Z(G)$ is an involution, then  $[i,i^g]=1$ since $z \in Z(G)$. Hence, $\zeta^i_0(g) = ii^g$. Since $\zeta_0^i(g)= iz$ by assumption, we have $i^g = z$. Thus $i = z$ and $\zeta^i_0(g) = 1$, a contradiction. \hfill $\Box$

\begin{lemma}\label{direct-heart}
Let $G=G_1 \times G_2$ be a direct product of finite groups $G_1$ and $G_2$. If $i=(i_1,i_2) \in G$ is an involution, then $\heart_i(G) = \heart_{i_1}(G_1) \times
\heart_{i_2}(G_2)$.
\end{lemma}
\proof Recall that the image of $\zeta^i_0$ belongs to the conjugacy classes of involutions in $C_G(i)$ by Theorem \ref{dist}. Since conjugacy classes of involutions in $G$ are direct products of conjugacy classes of $G_1$ and $G_2$,   the result follows.\hfill $\Box$

Let $j \in G$ be a pseudo-involution. Then we define $$\zeta^j_0(g) = {\rm i}(jj^g),$$ where $g \in G$ and ${\rm i}(jj^g)$ is an involution produced from $jj^g$ as in Section \ref{const-bb}. Observe that $\zeta^j_0(g) \in C_G(j)$. Moreover, we define $\heart_j(G)$ similarly for a pseudo-involution $j \in G$.

\begin{lemma}\label{heart-psl}
Let $G\cong \sl_2(q)$ and $j \in G$ be a pseudo-involution. Then $\heart_j(G)=Z(G)$.
\end{lemma}
\proof The result follows from the observation that $G$ has a unique central involution and the image of the function $\zeta^j_0$ is a set of involutions in $G$. \hfill $\Box$

\begin{lemma}\label{pseudo-sl2}
Let $G\cong \sl_2(q)$ and $j \in G$ be a pseudo-involution. Then $\langle \zeta_1^j(G)\rangle = N_G(\langle j \rangle)$ and $\langle \zeta_1^j(G)\rangle^{\prime \prime} =1$.
\end{lemma}
\proof We recall that if $jj^g$ has odd order $m$, then $\zeta_1^j(g) = (jj^g)^{(m+1)/2}g^{-1}$. Let $h=(jj^g)^{(m+1)/2}$, then it is straightforward to check that $j^h = j^2j^g$ which implies that $j^{hg^{-1}}=j^3$ since $j^2 \in Z(G)$. Hence, $\zeta_1^j(g) \in N_G(\langle j \rangle)$. Observe that if $j \in T$, then $N_G(\langle j\rangle) = N_G(T)$, $C_G(j) = C_G(T)=T$ and $|N_G(T)/C_G(T)|=2$. 
Hence, by Theorem \ref{dist}, we have $\langle \zeta_1^j(G)\rangle = N_G(\langle j \rangle)$.  Moreover, $\langle \zeta_1^j(G)\rangle^{\prime} \leqslant C_G(j)$ and $\langle \zeta_1^j(G)\rangle^{\prime \prime} =1$ since $C_G(j)$ is a cyclic group of order $q-1$ or $q+1$. \hfill $\Box$

\begin{lemma}\label{heart-symplectic}
Let $G\cong \ppsp_{2n}(q)$, $q>3$ and $i$ be an involution of type
$t_1$.
\begin{itemize}
\item[{\rm (a)}] If $n \geqslant 3$, then $\heart_i(G)^\prime \cong \sp_{2n-2}(q)$.
\item[{\rm (b)}] If $G\cong  \sp_4(q)$, then $\heart_i(G) = Z(G)$.
\item[{\rm (c)}] If $G\cong  \psp_4(q)$, then $\heart_i(G) \geqslant E(C_G(i))$, where $E(C_G(i))$ is the semisimple socle of $C_G(i)$.
\end{itemize}
\end{lemma}
\proof 
\begin{itemize}
\item[{\rm (a)}] Assume first that $G\cong  \sp_{2n}(q)$. Then $C_G(i)\cong \sl_2(q) \times \sp_{2n-2}(q)$. Let $V=V_- \oplus V_+$ be the decomposition of the corresponding vector space $V$, where $V_\pm$ are the eigenspaces of $i$ for the eigenvalues $\pm 1$. Then the dimension of $V_-$ is $2$ or $2n-2$. We assume that $\rm{dim} V_- =2$, the other case is analogous. It is clear that the dimension of the eigenspace for the eigenvalue $-1$ of the involution $\zeta^i_0(g)$ is at most 4 for any $g \in G$. Therefore the image of $\zeta^i_0$ contains non-central involutions in $\sp_{2n-2}(q)$ since $n \geqslant 3$. Hence, $\heart_i(G) \cong \{ \pm I_{2n} \} \times \sp_{2n-2}(q)$, where $I_{2n}$ is a $2n\times 2n$ identity matrix. In the case of $\psp_{2n}(q)$, we have $C_G(i) \cong \sl_2(q) \circ \sp_{2n-2}(q)$, and, by the same argument, all the involutions in the image $\zeta^i_0$ centralize the component isomorphic to $\sl_2(q)$ and do not centralize the other component since $n \geqslant 3$.

\item[{\rm (b)}] We have $C_G(i) \cong \sl_2(q) \times \sl_2(q)$. Hence, the set of involutions in $C_G(i)$ is $\{i,iz,z\}$, where $z\in Z(G)$ is the unique involution in $Z(G)$. By Lemma \ref{zeta0-lemma}, the image of $\zeta^i_0$ does not contain involutions $i$ and $iz$, and the result follows.

\item[{\rm (c)}] If $G\cong  \psp_4(q)$, then $C_G(i) \cong (\sl_2(q) \circ_2 \sl_2(q))\rtimes \langle t \rangle$, where $t$ is an involution interchanging the components. Since $\zeta^i_0$ does not produce the involution $i$ by Lemma \ref{zeta0-lemma}, it produces either an involution acting non-trivially on both components of $C_G(i)$ or the involution $t$. In either case, as $\heart_i(G)$ is normal subgroup in $C_G(i)$ by Theorem \ref{dist} (2), we conclude that $\heart_i(G) \geqslant E(C_G(i)) \cong \sl_2(q) \circ_2 \sl_2(q)$.\hfill $\Box$
\end{itemize}

The following Lemma is a direct consequence of Glauberman $Z^*$\/-\/Theorem \cite[page 262]{aschbacher2000} and will be used to prove the next theorem. 

\begin{lemma}\label{Glauberman}
Let $G$ be a non-abelian finite simple group and $i$ be an involution in $G$. Then there exists an involution $j \in C_G(i)$ such that $j\neq i$ and $j^g=i$ for some $g \in G$.
\end{lemma}

\begin{theorem}\label{heart}
Let $G$ be a finite simple group of Lie type over a field of odd characteristic $p$ and $i\in G$ be an involution. 
\begin{enumerate}
\item If $G$ is classical, then $\heart_i(G)$ contains the semisimple socle of $C_G(i)$ except when $G\cong \psp_{2n}(q)$ and $i$ is an involution of type $t_1$. 
\item If $G$ is exceptional, then $\heart_i(G)$ contains at least one component of $C_G(i)$.
\end{enumerate}
\end{theorem}
\proof  Let $i\in G$ be an involution. We prove the claim by constructing an involution in the image of $\zeta^i_0$ which does not centralize the component(s) in $C_G(i)$. It is clear that the existence of such involutions guarantees that $\heart_i(G)$ contains that the component(s) of the semisimple socle of $C_G(i)$ since $\heart_i(G)$ is a normal subgroup of $C_G(i)$ by Theorem \ref{dist}. 

Let $G\cong \psl_{n+1}(q)$, $n \geqslant 2$. Assume that $i$ is an involution of type $t_k$, where $2 \leqslant k \leqslant n/2$.  Then the semisimple socle $H$ of $C_G(i)$ is $H=H_1H_2$, where $H_1 \cong \sl_k(q)$ and $H_2 \cong \sl_{n+1-k}(q)$, and $[H_1,H_2]=1$. Observe that there is an involution $j \in C_G(i)$ of type $t_k$ in $G$ which acts as an involution of type $t_1$ in $H_1$ and of type  $t_{k-1}$ in $H_2$. Hence, $j=i^g$ for some $g \in G$. By construction, the involution $\zeta^i_0(g)=ii^g$ does not centralize the components $H_1$ and $H_2$. Hence, $H \leqslant \heart_i(G)$. 
Assume now that $i$ is an involution of type $t_1$ or $t_{(n+1)/2}^\prime$, then the semisimple socle of $C_G(i)$ is isomorphic to $\frac{1}{(n,q-1)} \sl_n(q)$ or  $\frac{1}{((n+1)/2,q-1)}\sl_{(n+1)/2}(q^2)$, respectively. In both cases, the involution $i$ is the unique involution in $Z(C_G(i))$, and, by Lemma \ref{Glauberman}, there is an involution $j=i^g \in C_G(i)$ for some $g \in G$, which does not centralize the semisimple socle of $C_G(i)$. Hence, the involution $\zeta^i_0(g) = i i^g$ is not central in $C_G(i)$ and the result follows. The proof for $\psu_n(q)$ is analogous.

Let $G\cong  \psp_{2n}(q)$ and $n \geqslant 3$. We refer to Lemma \ref{heart-symplectic}(c) for the case $n=2$, and Lemma \ref{heart-symplectic}(a) for the involutions of type $t_1$. Assume that $i \in G$ is an involution of type $t_k$, where $2\leqslant k \leqslant n/2$. Then the semisimple socle $H$ of $C_G(i)$ is $H = H_1 H_2$, where $H_1 \cong \sp_{2k}(q)$ and $H_2 \cong \sp_{2n-2k}(q)$. Take an involution $s_1 \in H_1$ which is of type $t_{k-1}$ in $G$ and an involution $s_2 \in H_2$ which is of type $t_1$ in $G$. It is clear that the involution $s=s_1s_2$ does not centralize the components $H_1$ and $H_2$. Moreover, we can assume that $s$ is of type $t_k$ in $G$ since $s_1$ and $s_2$ commute. Hence, $s = i^g$ for some $g \in G$. Now $\zeta^i_0(g) = ii^g \in H$ and $H \leqslant \heart_i(G)$.
If $i$ is an involution of type $t_n$ or $t_n^\prime$, then the semisimple socle of $C_G(i)$ is isomorphic to $\frac{1}{(2,n)}\sl_n^\eps(q)$, where $q \equiv \eps \, {\rm mod} \, 4$. In either case, that is for $\eps =\pm$, there exists $g \in G$ such that $i^g \in C_G(i)$  by Lemma \ref{Glauberman}. Note that $i^g$ is non-central in $C_G(i)$ as $i$ is the only involution in $Z(C_G(i))$. Hence, the result follows from a similar argument above.
   
Let $G\cong  \psizo_{2n+1}(q)$ and $n \geqslant 3$. For the proof of $n=2$, we refer to Lemma \ref{heart-symplectic} (c) as $\psizo_5(q) \cong \psp_4(q)$. Let $V$ be the natural module for $G$. Then the dimension of the eigenspace of each of the involutions $t_k$ or $t_k^\prime$, $1\leqslant k \leqslant n$, for the eigenvalue $-1$ is $2k$. We denote these involutions by $t_k$ to simplify the notation, see Table \ref{cent-table} for the structure of their centralizers. 
Let $i$ be an involution of type $t_1$, then the semisimple socle $H$ of $C_G(i)$ is isomorphic to $H\cong \psizo_{2n-1}(q)$. Notice that there is a non-central involution of type $t_1$ in $H$, which is conjugate to $i$ in $G$. Hence, $H \leqslant \heart_G(i)$ by similar arguments above.
Assume now that $i$ is an involution of type $t_k$, $2 \leqslant k < n$. Then the semisimple socle of $C_G(i)$ is $H \cong \psizo^\varepsilon_{2k}(q) \times \psizo_{2(n-k)+1}(q)$, where $q^k \equiv \varepsilon \, {\rm mod} \,4$. Note that $\psizo^\varepsilon_{2k}(q)$ contains an involution $s_1$, which is conjugate to an involution of type $t_{k-1}$ in $G$, and $\psizo_{2(n-k)+1}(q)$ contains an involution $s_2$, which is conjugate to an involution of type $t_1$ in $G$. Since $s_1$ and $s_2$ commute, we can assume that the involution $s = s_1 s_2$ is conjugate to an involution of type $t_k$ in $G$, and $s = i^g$ for some $g \in G$. Now $\zeta^i_0(g) = ii^g$ does not centralize the components of $H$ since $s_1$ and $s_2$ are non-central involutions in the corresponding components and the result follows.
If $i\in G$ is an involution of type $t_n$, then, following the notation of the proof of Lemma \ref{heart-symplectic}(a), we have ${\rm dim}V_{-}=2n$. The semisimple socle $H$ of $C_G(i)$ is isomorphic to $ \psizo^\varepsilon_{2n}(q)$, where $q^n \equiv \varepsilon \, {\rm mod} \,4$. Observe that $ii^g \in C = C_{\sl(V)}(V_{-} \cap V_{-}g)$ for any $g \in G$. Since ${\rm dim}V_{-} =2n$, we have that $V_{-} \cap V_{-}g$ has codimension $\leqslant 2$ and $C/O_p(C)$ is a subgroup of $\sl_2(q)$. Notice also that we can choose $g \in G$ so that $i^g \in C_G(i)$ by Lemma \ref{Glauberman}. Hence, $\zeta^i_0(g) = ii^g$ and it is clear that $\zeta^i_0(g)$ does not centralize $H$.

Let $G\cong \pso_{2n}^\eps(q)$ and $n \geqslant 4$. For the involutions of types $t_k$ or $t_k^\prime$, where $1\leqslant k <n-1$, the proof is similar to the case $G\cong \psizo_{2n+1}(q)$. If $i$ is of type $t_{n-1}$ or $t_n$, then the semisimple socle of $C_G(i)$ is isomorphic to $\frac{1}{2}\sl_n^\eps(q)$, where $q \equiv \eps \, {\rm mod} \, 4$. In either case, by Lemma \ref{Glauberman}, there exists $g \in G$ such that $i^g \in C_G(i)$ and $i^g$ is non-central in $C_G(i)$ since $i$ is the only involution in $Z(C_G(i))$.

Let $G\cong \/^3D_4(q)$ or $G_2(q)$, then there is only one conjugacy class of involutions in $G$, and $C_G(i) \cong \sl_2(q) \circ_2 \sl_2(q^3)$ or $ \sl_2(q) \circ_2 \sl_2(q)$, respectively. Notice that there is an involution $s \in C_G(i)$ which does not centralize the components, and it is necessarily conjugate to $i$ since there is only one conjugacy class of involutions. Let $s=i^g$ for some $g \in G$, then $\zeta^i_0(g)=ii^g$ does not centralize the components, hence the result follows. We apply the same argument to the group $^2G_2(q)$ as there is only one conjugacy class of involutions in $^2G_2(q)$.

Let $G$ be a group of type $F_4(q)$. Let $i$ be an involution of type $t_4$, then the semisimple socle $H$ of $C_G(i)$ is isomorphic to  $H \cong {\rm Spin}_9(q)$ and $Z(C_G(i)) = \langle i \rangle$. By Lemma \ref{Glauberman}, there exists $g \in G$ such that $i\neq i^g \in C_G(i)$ and $\zeta^i_0(g) = i i^g$ is non-central in $H$. If $i$ is an involution of type $t_1$, then $C_G(i) \cong \sl_2(q) \circ_2 \sp_6(q)$ and $Z(C_G(i)) = \langle i \rangle$. Again, by Lemma \ref{Glauberman}, there exists $g \in G$ such that $\zeta^i_0(g) = i i^g$ does not centralize $\sp_6(q)$. Hence, $\sp_6(q) \leqslant \heart_i(G)$. By the same argument, if $\sl_2(q)$ is a component in the semisimple socle of $C_G(i)$ for $G\cong E_6(q)$, $E_7(q)$, $E_8(q)$, then $\heart_i(G)$ contains the other (quasi)simple component of the semisimple socle of $C_G(i)$. 

If the centralizers of involutions in the remaining exceptional groups $E^{\varepsilon}_6(q)$, $E_7(q)$ and $E_8(q)$ do not have a component isomorphic to $\sl_2(q)$, then the semisimple socle is a (quasi)simple group, see Table \ref{cent-table}. Hence, we can apply previous arguments to these cases.\hfill $\Box$

\begin{corollary}\label{heartcor}
Let $G$ be a quasi-simple group of Lie type defined over a field of odd characteristic and $i \in G$ be a non-central involution in $G$. 
\begin{enumerate}
\item If $G$ is classical, then $\heart_i(G)$ contains the semisimple socle of $C_G(i)$ except when $G\cong \ppsp_{2n}(q)$ and $i$ is an involution of type $t_1$. 
\item If $G$ is exceptional, then $\heart_i(G)$ contains at least one component of $C_G(i)$.
\end{enumerate}
\end{corollary}
\proof  Note that the non-central involutions in $G$ map to involutions in $G/Z(G)$ and we apply the same arguments as in Theorem \ref{heart}. \hfill $\Box$

\begin{remark} {\rm In our algorithms, the map $\zeta_0^i$ will be considered only in the construction of at least one quasisimple component of a centralizer of an involution. Therefore, the exceptions in Theorem \ref{heart} (2) and Corollary \ref{heartcor} (2) do not provide any difficulties in our main algorithm. }
\end{remark}

\section{Construction of a long root $\sl_2(q)$-subgroup}\label{const-lr-sl2}
The aim of this section is to prove Theorem \ref{main1}. %All algorithms in this section except Algorithm \ref{g2and3d4} are outlined in \cite{suko01}. Here we present them in detail.
We present the following algorithm.

\begin{algo}\label{main-simple-alg} {\rm (cf. \cite[Algorithm 5.1]{suko01})}\textsc{ ``Construction of a long root $\sl_2(q)$-subgroup in a finite simple group of Lie type''}
\begin{itemize}
\item[] Input: A black-box group isomorphic to a finite simple group $G$ of Lie type defined over a field of odd order $q>3$ except $\psl_2(q)$ and $^2G_2(q)$. 
\item[] Output: A black-box group which is a long root $\sl_2(q)$-subgroup in $G$.
\end{itemize}
\end{algo}

We need the following terminology. Let $G$ be a group and $G_i \leqslant G$, $i=1,\ldots ,n$. Assume that 
$$G= \langle G_i \mid i=1,\ldots, n \rangle$$
and $G_k$ commutes with $G_l$ for any $k\neq l$. Then we say that $G$ is a \textit{commuting product} of the subgroups $G_i$ for $i=1,\ldots ,n$. 

Algorithm \ref{main-simple-alg} has three major components:
\begin{itemize}
\item[1.] Construct the centralizers of involutions recursively to find a commuting product of the groups $\ppsl_2(q^k)$ in $G$, where $k\geqslant 1$ may vary; Section \ref{find-prosl2}.
\item[2.] Construct one of the components $K\cong \ppsl_2(q^k)$ in the commuting product found in Step 1; Section \ref{find-sl2}.
\item[3.] Check if $K$ is a long root $\sl_2(q)$-subgroup of $G$; Section \ref{long-root-sl2-subgroup}.
\end{itemize}

\begin{remark}
\noindent {\rm {\bf 1.} In Step 1, we may have a commuting product of the subgroups $\ppsl_2(q)$, where $q$ may vary but the characteristic of the underlying field is the same, for example, in $G\/=\/^3D_4(q)$ there is only one class of involutions and $C_G(i) \cong \sl_2(q) \circ_2 \sl_2(q^3)$ for an involution $i\in G$.

\noindent {\bf 2.} The long root $\sl_2(q)$-subgroups in simple groups of Lie type of odd characteristic are indeed isomorphic to $\sl_2(q)$ by Theorem \ref{long-root-sl2}. Therefore if we obtain a direct product of subgroups $\psl_2(q)$ in Step 1, then we conclude that it does not contain a long root $\sl_2(q)$-subgroup as a component and we repeat Step 1 to construct a new commuting product of subgroups $\ppsl_2(q)$.

%\noindent{\bf 3.} If $G\cong G_2(q)$ or $^3D_4(q)$, then we have a special subroutine to construct a long root $\sl_2(q)$-subgroup, see Subsection \ref{sec-g2and3d4}.  
%\noindent{\bf 3.} If we have a subgroup $K$ which is not a long root $\sl_2(q)$-subgroup in $G$, then we construct an another component in the commuting product constructed in Step 3. 
}
\end{remark}

It is easy to see that the proof of Theorem \ref{main1} follows from Algorithm \ref{main-simple-alg} and  Corollary \ref{semisocle2}.

\subsection{Constructing products of $\ppsl_2(q)$}\label{find-prosl2}
In this section we present an algorithm concerning Step 1 of Algorithm \ref{main-simple-alg}. The algorithm is as follows.
%Our aim in this section is to construct commuting product of subgroups $\ppsl_2(q)$. 
%We achieve this by constructing centralizers of involutions recursively. For that we use both of the functions $\zeta_0^i$ and $\zeta_1^i$ for the involutions $i \in G$.

\begin{algo}\label{prod-sl2q}{\rm (cf. \cite[Algorithm 5.3]{suko01})} \textsc{``Construction of a commuting product of subgroups $\ppsl_2(q)$''}
\begin{itemize}
\item[] Input: A black-box group isomorphic to a finite simple group $G$ of Lie type defined over a field of odd order $q>3$ except $\psl_2(q)$ and $^2G_2(q)$.
\item[] Output: A black-box group which is a commuting product of subgroups $\ppsl_2(q^k)$ for various $k\geqslant 1$.
\end{itemize}
\end{algo}

\noindent\textsc{Description of the algorithm:}
\begin{itemize}
\item[1.] Produce an involution $i={\rm i}(g)$ from a random element $g\in G$.
\begin{itemize}
\item[$\bullet$] Check if $i\in Z(G)$. If we always have $i \in Z(G)$ after $O(N)$ steps, where $N$ is the input length of the black-box group $G$, then return $G$. Otherwise, go to the next step.
\end{itemize}
\item[2.] Construct a subgroup $L$ where $\heart_i(G) \leqslant L \leqslant C_G(i)$ by using $\zeta^i=\zeta_0^i\sqcup \zeta_1^i$.
\item[3.] Construct $L^\prime$ and $L^{\prime \prime}$.
\begin{itemize}
\item[$\bullet$] If $L^{\prime \prime} = 1$ then return $\langle i^G\rangle^{\prime}$. Otherwise, set $G=L^{\prime \prime}$ and go to Step 1.
\end{itemize}
\end{itemize}

\subsubsection{The involution $i={\rm{i}}(g)$ from a random element $g$}
We pick a random element $g \in G$ and check if $g^{2^k}=1$, where $E=2^km$, $m$ odd, is the exponent of the group given as an input.  By \cite[Corollary~5.3]{isaacs95.139}, the proportion of elements of even order in groups of Lie type of odd characteristic is at least 1/4. As soon as we find an  element $g \in G$ of even order, we construct $i={\rm i}(g)$ as described in Section \ref{const-bb}.

Next, we need to check whether $i={\rm{i}}(g)$ is central in $G$ or not. Notice that this procedure is unavoidable in the recursive steps of our algorithm as $C_G(i)^{\prime \prime}$ contains central involutions in most of the cases, see Table \ref{cent-table}. Since the generators of $G$ are given, it is straightforward to decide whether $i\in Z(G)$ or not. We can find a non-central involution in view of the following lemma.

\begin{lemma}\label{cent-invo}
Let $G$ be a universal version of a finite group of Lie type defined over a field of odd characteristic and $G\ncong \sl_2(q)$. Then the proportion of elements in $G$ producing non-central involutions is bounded from below by a function of the Lie rank of $G$.
\end{lemma}
\proof Let $T$ be a maximal torus twisted by an element $w \in W$ as in Section \ref{maximal_tori}. Then the probability of being an element in a torus conjugate to $T$ depends only on the Weyl group $W$ not on the order of the field by the remark after Equation (\ref{toriat-q}) in Section \ref{maximal_tori}. Therefore if $T$ is any torus and $i \in T$ an involution, then the number of elements $g\in G$ with $i={\rm i}(g)$ depends only on $W$ not on the order of the field. In classical groups, tori twisted by the cycles of length less than the length of the longest cycle do not contain central involutions. This observation immediately follows from the tables of the centralizers of involutions and the orders of the corresponding tori given in \cite{carter70.g1}. Note that among the exceptional groups only the universal version of $E_7(q)$ has a central involution in which case the same arguments apply. \hfill $\Box$

We give an example in the easiest case. Let $G\cong \sl_n(q)$, $n$ even and $n\geqslant 5$. Then a maximal torus $T$ twisted by a product of an $(n-2)$-cycle and a $2$-cycle is of the form $T=\frac{1}{q-1}(T_1 \times T_2)$, where $T_1$ is a cyclic group of order $(q^{n-2}-1)$ and $T_2$ is a cyclic group of order $(q^2-1)$. Note that as $n$ is an even number, $(q^2-1)$ divides $(q^{n-2}-1)$ and therefore involutions produced from random elements in $T$ belong to $T_1$ with probability very close to 1. It can be observed from Table \ref{cent-table} that $T_1$ has involution $i$ of type $t_2$ in $G$ by comparing the orders of $C_G(i)$ and $T_1$. Since $W={\rm{Sym}}(n)$, the symmetric group on $n$ letters, and $n\geqslant 5$, $C_W(w)$ has order $2(n-2)$. Therefore the probability of producing non-central involutions from random elements is at least $\frac{1}{2(n-2)}$.

%Hence after $O(n)$ repetitions we can find a non-central involution $i={\rm{ i}}(g) \in G$ by Lemma \ref{cent-invo}. If we are always in a situation of producing central involutions then 

Thus, if we can not find a non-central involution of the form $i={\rm i}(g)\in G$ in $O(N)$ tests of random elements $g \in G$, then by \cite{griess78.241}, we can deduce that $G$ is isomorphic to the direct product of copies of $\sl_2(q)$, and we take this subgroup as an output. Otherwise, we go to the next step. Note that the number,  $O(N)=O(n^2\log q)$, of repetitions can be reduced to $O(n)$, if the Lie rank of $G$ is known.

\subsubsection{Construction of $C_G(i)$}\label{concgi}
Let $i={\rm{i}}(g)$ be a non-central involution in $G$. Take a subset $S \subset G$ consisting of random elements of $G$ and consider the subgroup $\langle \zeta^i(S)\rangle$. 

%\begin{theorem}{\rm (\cite{kantor90.67, liebeck95.103})}\label{prob-gen}
%The probability that two randomly chosen elements in a finite simple group $G$ generate $G$ tends to $1$ as $|G|$ tends to $\infty$.
%\end{theorem}

%The following result is the combination of \cite[Proposition 10]{kantor90.67} and  \cite{liebeck95.103}.

%\begin{theorem}\label{direct-prob-gen}
%Let $G$ be a direct product of finite simple groups of Lie type. Then the probability that two elements of $G$ generate $G$ approaches $1$ as the orders of each of the factors of $G$ approach $\infty$.
%\end{theorem}

By Lemma \ref{semisocle}, the semisimple socle of $C_G(i)$ is $C_G(i)^{\prime \prime}$, which is a commuting product of (quasi)simple groups of Lie type. By \cite{kantor90.67, liebeck95.103}, randomly chosen two elements in a finite simple group $G$ generate $G$ with probability tending to $1$ as $|G| \rightarrow \infty$. By \cite[Proposition 10]{kantor90.67}, a similar result holds for a direct product of simple groups where the orders of each of the factors approach $\infty$. Thus randomly chosen two elements in $C_G(i)$ generate a subgroup containing the semisimple socle of $C_G(i)$ with probability close to 1 when the size of the field is large. We know that the map $\zeta_1^i$ produces uniformly distributed random elements in $C_G(i)$ by Theorem \ref{dist}. Therefore if a reasonable number of elements $g\in S$ satisfy $\zeta_1^i(g) \neq 1$, then we conclude that $\langle \zeta_1^i(S) \rangle$ contains the semisimple socle of $C_G(i)$.

If $\zeta^i_1$ is not defined for the elements in $S$, then $\zeta_0^i$ is defined for all the elements in $S$. In this case, we use only the function $\zeta^i_0$ for the generation of a subgroup in $C_G(i)$. Recall that the image of the map $\zeta_0^i$ is a normal subset of $C_G(i)$ by Theorem \ref{dist}. By \cite[Theorem 1.1]{liebeck95.103}, for any normal subset $S\subset G$ of a finite simple group $G$, there exists a constant $c$ such that $S^m=G$ for any $m\geqslant c  {\rm{log}}|G|/{\rm{log}}|S|$. A similar result also holds for a direct product of simple groups by considering direct product of normal subsets. Therefore, in such a case, we take $S$ to be sufficiently large so that we have $\langle \zeta_0^i(S)\rangle = \heart_i(G)$. Now, by Lemma \ref{direct-heart} and Theorem \ref{heart}, $\heart_i(G)$ contains at least one component of $C_G(i)$. Experiments in GAP show that for a reasonably sized subset $S\subset G$ we have $\langle \zeta^i_0(S)\rangle = \heart_i(G)$. A reasonable number of generators for $\heart_i(G)$ produced by $\zeta_0^i$ is, for example, $50$ for groups of $50 \times 50$ matrices. %Recall that the image of the map $\zeta_0^i$ is a normal subset of $C_G(i)$ by Theorem \ref{dist}. Here we note two results of Liebeck and Shalev on the generation of finite simple groups by conjugacy classes. 

%\begin{theorem}{\rm (\cite{liebeck01.383})}
%Let $G$ be a finite simple group and $S \subseteq G$ a normal subset. Then there exists a constant $c$ such that $S^n = G$ for any $n\geqslant c  {\rm{log}}|G|/{\rm{log}}|S|$.
%\end{theorem}

%\begin{theorem}{\rm (\cite{liebeck01.383})}
%There exists a constant $c$ such that any element of any finite simple group $G$ can be written as a product of $c$ involutions. 
%\end{theorem}

\subsubsection{Semisimple socle of $C_G(i)$}

From now on we assume that  $\heart_i(G) \leqslant \langle \zeta^i(S)\rangle$, and by definition $\langle \zeta^i(S)\rangle \leqslant C_G(i)$.

We shall construct the semisimple socle of $\langle \zeta^i(S)\rangle$. Note that the derived subgroup of a black-box group can be constructed in Monte--Carlo polynomial time \cite[Corollary 1.6]{babai95.296}. By Lemma \ref{semisocle} and Theorem \ref{heart}, the subgroup $\langle \zeta^i(S)\rangle^{\prime \prime}$ contains at least one of the components of $C_G(i)^{\prime \prime}$ unless  $G\cong  \sp_4(q)$, $i\in G$ is an involution of type $t_1$ and $\langle \zeta^i(S)\rangle^{\prime \prime} = \heart_i(G)^{\prime \prime}$, see Lemma \ref{heart-symplectic}. However, the probability that $\langle \zeta^i(S)\rangle^{\prime \prime} = \heart_i(G)^{\prime \prime}$ in $G\cong \sp_4(q)$ for a classical involution $i$ is close to $0$ by Theorem \ref{zeta1-lemma}. 

Assume that we are at the $k^{\rm th}$  recursion of our algorithm, namely 
$$C_{i_k} = \langle\zeta^{i_k}(S_k)\rangle^{\prime \prime},$$
where $i_k={\rm{i}}(g)$, $g \in C_{i_{k-1}}^{\prime \prime}$ and $S_k \subset C_{i_{k-1}}^{\prime \prime}$. Set 
$H = C_{i_{k-1}}^{\prime \prime}.$

If $C_{i_k}^{\prime \prime} \neq 1$, then we set $G=  C_{i_k}^{\prime \prime}$ and go to Step 1 of Algorithm \ref{prod-sl2q}. If $C_{i_k}^{\prime \prime} =1$, then by Lemmas \ref{direct-heart}, \ref{heart-psl}  and \ref{pseudo-sl2}, $H$ is a commuting product of subgroups isomorphic to $\ppsl_2(q)$ and $i_k$ acts as a pseudo-involution on some components of $H$. Hence, we return  the subgroup $\langle i_k^H \rangle^\prime$. Note that the normal closure of a black-box group can be constructed in Monte--Carlo polynomial time \cite[Theorem 1.5]{babai95.296}.

\subsection{Constructing $\sl_2(q)$}\label{find-sl2}
The aim of this section is to present an algorithm concerning Step 2 of Algorithm \ref{main-simple-alg}. We construct a normal subgroup $\sl_2(q^k)$ which appears as a component in a given commuting product of groups $\ppsl_2(q^l)$ for various $l$ found in Section \S \ref{find-prosl2}.

We need the following results.

\begin{lemma}\label{pseudo-prod-sl2}
Let $K=\sl_2(q)$, $q>3$, and $L=K\times K$. Then the probability of producing a pseudo-involution acting trivially on one component is at least $1/24$.
\end{lemma}

\proof  It is well known that all semisimple elements in $K \backslash Z(K)$ are regular, therefore they belong to only one torus. There are two conjugacy classes of tori in $K$: split and non-split tori. A split torus is a cyclic group of order $q-1$ and a non-split torus is a cyclic group of order $q+1$. 

Let $T_1$ be a split torus in $K$, then 
$$|K:N_K(T_1)| = q(q+1)/2.$$ 
Let $T_2$ be a non-split torus in $K$, then 
$$|K:N_K(T_2)| = q(q-1)/2.$$ 
Hence, the total number of tori is $q^2$. Therefore the probability of a semisimple element belonging to a split torus  is 
$$\frac {1}{2} \frac{q(q+1)}{q^2} > \frac{1}{2} $$
and to a non-split torus  is
$$\frac {1}{2} \frac{q(q-1)}{q^2} \approx \frac{1}{2} > \frac{1}{3}$$
since $q>3$. Let $g=(g_1,g_2) \in L$. Then $g_1^q$ and $g_2^q$ belong to different classes of tori in $L$ with probability at least $1/6$. Note that the order of one of the tori is divisible by $4$ and the probability of finding an element in this torus which has order divisible by $4$ is at least $1/4$. Therefore the probability that a pseudo-involution produced from a random element acts non-trivially on only one component of $L$ is at least $1/24$. \hfill $\Box$

\begin{corollary}\label{cor-pseudo-prod-sl2}
Let $L$ be a commuting product of the groups isomorphic to $\sl_2(q)$. Then the probability of producing a pseudo-involution acting trivially on some components of $L$ is at least $1/24$.
%we can find a pseudo-involution which acts trivially on some components of $L$ with probability at least $1/24$.
\end{corollary}

\begin{lemma}\label{pseudo-zeta1}
Let $K=\sl_2(q)$ and $t \in K$ be a pseudo-involution. Then elements of the form $tt^g$ have odd order with probability at least $1/6$.
\end{lemma}

\proof The elements $z=tt^g$ belongs to a torus of order $q-1$ or $q+1$.  The probability that $z$ belongs to a certain type of torus is at least $1/3$ by the proof of Lemma \ref{pseudo-prod-sl2}, and one of $(q\pm 1)/2$ is odd. Therefore the element $z$ has odd order with probability at least $1/6$. \hfill $\Box$

Now we can present our algorithm.

\begin{algo}\label{sl2q}{\rm (cf. \cite[Algorithm 5.4]{suko01})} \textsc{``Construction of $\sl_2(q)$''}
\begin{itemize}
\item[] Input: A black-box group $L$ which is a commuting product of groups $\ppsl_2(q^l)$ for various $l$.
\item[] Output: A black-box group $\sl_2(q^k)$ for some $k$ appearing as a component in the commuting product or return the statement ``$L$ is a direct product of groups $\psl_2(q^l)$ for various $l$''.
\end{itemize}
\end{algo}

\noindent\textsc{Description of the algorithm:}
\begin{itemize}
\item[1.] Produce a pseudo-involution $t \in L$. If we can not find any pseudo-involution then return the statement ``$L$ is a direct product of $\psl_2(q^l)$ for various $l$'' and go to the beginning of Step 1 of Algorithm \ref{main-simple-alg}. Otherwise, go to the next step.
\item[2.] Construct $\langle t^L \rangle^\prime$ and check if $\langle \zeta_1^t(S)\rangle^{\prime \prime} =1$ for a reasonably sized subset $S\subset L$. If $\langle \zeta_1^t(S)\rangle^{\prime \prime} \neq 1$, then set $L=\langle t^L \rangle^\prime$ and go to Step 1.
\end{itemize}

{\bf Step1:} Recall that the random elements of $\sl_2(q)$ belong to a torus of order $q \pm 1$ with probability $1-O(1/q)$ by \cite{guralnick01.169}. Observe that one of the numbers $(q\pm1)/2$ is even and at least half of the elements in such  tori have order a multiple of 4. Therefore the probability of finding an element of order a multiple of $4$ is close to $1/4$. Hence, if $L$ has a component isomorphic to $\sl_2(q)$, then we can produce a pseudo-involution from random elements by following the procedure in Section \ref{const-bb}. If we can not find a pseudo-involution in $L$, then we deduce that $L$ is a direct product of $\psl_2(q)$ and we start the procedure from the beginning.

{\bf Step2:} Let $t$ be a pseudo-involution produced from a random element in $L$. Then
\[
L=\langle t^L \rangle^\prime C_L(t)^{\prime \prime}.
\]
By Corollary \ref{cor-pseudo-prod-sl2}, $t$ acts trivially on some components of $L$ with probability at least $1/24$. Hence, the number of components in the subgroup $\langle t^L \rangle^\prime$ is less than the number of components of $L$.  We know that $ \langle \zeta_1^t(S)\rangle^{\prime \prime}  \leqslant C_L(t)^{\prime \prime} $ for any subset $S\subset L$ by Lemma \ref{pseudo-sl2}.  By using Lemma \ref{pseudo-zeta1}, we can choose a reasonably sized subset $S\subset L$ so that $C_L(t) = \langle \zeta_1^t(S)\rangle$.  Now if $\langle \zeta_1^t(S)\rangle^{\prime \prime} \neq 1$, then we set $L=\langle t^L \rangle^\prime$  and continue in this way. Note that if $\langle \zeta_1^t(S) \rangle^{\prime \prime}=1$, it may not necessarily be true that $L\cong \sl_2(q^k)$ since $t\in L$ may act non-trivially on more than one component of $L$. Therefore we check if $\langle \zeta_1^t(S) \rangle^{\prime \prime}=1$ for different pseudo-involutions $t \in L$, say $m$ times. If we always have $\langle \zeta_1^t(S) \rangle^{\prime \prime}=1$ for different pseudo-involutions $t \in L$, then we deduce that $L\cong \sl_2(q^k)$ for some $k$ in the commuting product where the probability of error is at most $(1-1/24)^m$  by Corollary \ref{cor-pseudo-prod-sl2}. 

\begin{remark}\label{comm-sl2}
{\rm  Notice that one can extend Algorithm \ref{sl2q} to construct each component isomorphic to $\sl_2(q)$ in a commuting product $L$ of subgroups $\ppsl_2(q)$. For this, let $K\cong \sl_2(q)$ be a component of $L$ and $t \in K$ a pseudo-involution. Then we set $L=C_G(t)^{\prime \prime}$ and apply Algorithm \ref{sl2q}. }
\end{remark}

%\newpage
\subsection{A long root $\sl_2(q)$}\label{long-root-sl2-subgroup}
%The aim of this section is to prove Theorem \ref{main2} by exhibiting the algorithm. 

%\begin{algo}\label{long-sl2q}{\rm (cf. \cite[Algorithm 5.6]{suko01})} \textsc{``Checking whether a given $\ppsl_2(q)$ is a long root $\sl_2(q)$''}
%\begin{itemize}
%\item[] Input: A black-box subgroup $K\leqslant G$ which is known to be isomorphic to $\sl_2(q)$.
%\item[] Output: The truth value of the statement: ``$K$ is a long root $\sl_2(q)$-subgroup in $G$''.
%\end{itemize}
%\end{algo}

%\noindent\textsc{Description of the algorithm:}\\
%\noindent\textbf{Step 1.} Construct an involution $z\in K$ and $C=C_G(z)^{\prime\prime}$.\\
%\noindent\textbf{Step 2.} Check whether $\langle K,K^g\rangle = K$ for random $g \in C$.

In this section we present an algorithm concerning Step 3 of Algorithm \ref{main-simple-alg}. We determine whether the subgroup $K\cong \sl_2(q^k)$ constructed in Section \ref{find-sl2} is a long root $\sl_2(q)$-subgroup in $G$ where $G \ncong\/ ^3D_4(q)$ or $G_2(q)$. The groups $G_2(q)$ and $^3D_4(q)$ are treated in Subsection \ref{sec-g2and3d4}.

\begin{algo}\label{long-sl2q}{\rm (cf. \cite[Algorithm 5.6]{suko01})} \textsc{``Checking whether a given $\ppsl_2(q)$ is a long root $\sl_2(q)$''}
\begin{itemize}
\item[] Input: A black-box subgroup $K\leqslant G$ which is known to be isomorphic to $\sl_2(q)$ where $G \ncong\/ ^3D_4(q)$ or $G_2(q)$.
\item[] Output: The truth value of the statement: ``$K$ is a long root $\sl_2(q)$-subgroup in $G$''.
\end{itemize}
\end{algo}

\noindent\textsc{Description of the algorithm:}
\begin{itemize}
\item[1.] Find the field size $q$.
\item[2.] Construct $C=C_G(z)^{\prime \prime}$, where $z\in Z(K)$. 
\item[3.] Construct $N=\langle K, K^g\rangle$ for a random $g \in C$.
\item[4.] Check if $n^{q(q^2-1)}\neq 1$ for random elements $n \in N$. If such an element is found, then return the statement ``$K$ is not a long root $\sl_2(q)$-subgroup''. If we can not find such an element then we deduce that $K$ is a long root $\sl_2(q)$-subgroup.
\end{itemize}

Recall that if $K$ is a long root $\sl_2(q)$-subgroup in $G$, then, by Theorem \ref{long-root-sl2}, $K = K^g$ for any $g \in C=C_G(z)^{\prime \prime}$, where $z\in Z(K)$. If not, then we prove that $K$ is strictly contained in the subgroup $N=\langle K, K^g \rangle$ for a random $g \in C$ with probability close to 1 except for the groups $G_2(q)$ and $^3\/D_4(q)$.

First, we find the size of the field by using \cite[Algorithm 5.5]{suko01}. Then we construct $C=C_G(z)^{\prime\prime}$ by using $\zeta_0^z \sqcup \zeta_1^z$ and consider $N = \langle K,K^g \rangle$ for random $g \in G$.

%We look for an element $n \in N$ such that  $n^{q(q^2-1)}\neq 1$. As soon as we find such an element, we conclude that $K$ is not a long root $\sl_2(q)$-subgroup. If we always have $n^{q(q^2-1)}=1$, then we conclude that $K$ is a long root $\sl_2(q)$-subgroup with probability close to 1. 

If $K$ is not a long root $\sl_2(q)$-subgroup then either $K$ is a short root $\sl_2(q)$-subgroup or the order of the field is increased in one of  the recursive steps of Algorithm \ref{prod-sl2q}. 

Assume first that the order of the field is increased in the recursive construction of centralizers of involutions. This is possible only when $G\cong \psl_n^\eps(q)$ where $n$ is an even integer, see Table \ref{cent-table}. Assume that $G\cong \psl_n(q)$ and $K \cong \sl_2(q^{2^k})$ for some $k\geqslant 1$. Then $K$ can be embedded naturally in $\sl_{2^{k+1}}(q)$ and $C \cong \sl_{2^{k+1}}(q) \circ \sl_{n-2^{k+1}}(q)$. Now take $K_1 \cong \sl_2(q^{2^{k-1}}) < K$. Then  $N_1=\langle K_1, K_1^g \rangle \cong \sl_4(q^{2^{k-1}})$ with probability at least $1-O(1/q^{2^{k-1}})$ for $g \in C$ by Lemma \ref{psl}. Now $N_1$ contains sufficiently many elements whose orders do not divide $|K|=q^{2^k}(q^{2^{k+1}}-1)$. The case $G\cong \psu_n(q)$ is analogous.

Now we assume that $K$ is a short root $\sl_2(q)$-subgroup in $G$. Then by \cite[Table 14.4]{aschbacher77.353}, $K\cong \psl_2(q^2)$, $\psl_2(q)$ or $\psl_2(q^2)$ for $G\cong  \psu_n(q)$, $\psizo_{2n+1}(q)$, $\pso_{2n}^-(q)$, respectively. Notice that it is easy to recognize these cases since $K$ does not contain a central involution. Now we are left with the cases $G\cong \psp_{2n}(q)$, $F_4(q)$ and $^2E_6(q)$.

If $G\cong  \psp_{2n}(q)$, then $K\cong\sl_2(q)$ and $C\cong \sp_4(q)\circ_2\sp_{2n-4}(q)$. Here $K$ is contained in the component isomorphic to $\sp_4(q)$. We have $N=\langle K, K^g \rangle \cong \sp_4(q)$ with probability at least $1 - 1/q$ for random $g \in C$, which follows from a similar idea of the proof of Lemma \ref{psp}. It is clear that the subgroup $N$ contains sufficiently many elements whose orders do not divide $q(q^2-1)$. 

If $G\cong F_4(q)$, then $K\cong \sl_2(q)$ and $C\cong  {\rm Spin}_9(q)$ \cite[29.7]{aschbacher77.353}. Since $C/Z(C)\cong \psizo_9(q)$, it is enough to obtain the estimates in the classical group $\psizo_9(q)$. Let $V$ be the natural module for $\psizo_9(q)$ and $K$ be a short root $\sl_2(q)$-subgroup in $\psizo_9(q)$, then $[K,V]$ is a orthogonal $3$-space of Witt index 1. Following the same idea in the proof of Lemma \ref{psizo}, we obtain that $N/Z(N) \cong \pso_6^+(q)$ with probability at least $1-O(1/q)$, where $N=\langle K,K^g \rangle$ for $g\in C$. Now the subgroup $N$ contains sufficiently many elements whose orders do not divide $q(q^2-1)$. 

If $G\cong \/^2\/E_6(q)$, then $K\cong \sl_2(q^2)$ and $C\cong {\rm Spin}_{10}^-(q)$ \cite[29.7]{aschbacher77.353}. Hence, $C/Z(C)\cong \pso_{10}^-(q)$. In this case it is enough to make the estimates in $\psizo_{10}^-(q)$. Notice that $[K,V]$ is an orthogonal $4$-space of Witt index 1, where $V$ is the natural module for $\psizo_{10}^-(q)$ and $K$ is a short root $\sl_2(q)$-subgroup in $\psizo_{10}^-(q)$. The rest is similar to the arguments above.

If we can not find an element $n \in N$ satisfying $n^{q(q^2-1)}\neq 1$, then we conclude that $K$ is a long root $\sl_2(q)$-subgroup of $G$. Notice that the output ``$K$ is not a long root $\sl_2(q)$-subgroup'' is always true.

%Note that if $G\cong \/^3D_4(q)$ or $G_2(q)$, then $G$  has only one conjugacy class of involutions and $C \cong \sl_2(q) \circ_2 \sl_2(q^3)$ or $\sl_2(q) \circ_2 \sl_2(q)$, respectively. In these cases we have $K = K^g$ for any $g \in C_G(z)^{\prime\prime}$. Thus this approach fails to recognize long root $\sl_2(q)$-subgroups in the groups $G_2(q)$ and $^3\/D_4(q)$, and we present a special routine in these cases in the next subsection. 

%If Algorithm \ref{long-sl2q} returns the statement ``$K$ is not a long root $\sl_2(q)$-subgroup'', then we construct another component in the commuting product $L$ of subgroups isomorphic to $\ppsl_2(q)$ constructed by Algorithm \ref{prod-sl2q} by following the procedure in Remark \ref{comm-sl2}. If $K=L$, then we start the procedure from the beginning. 

\subsubsection{The groups $G_2(q)$ and $^3\/D_4(q)$}\label{sec-g2and3d4}
Let $G\cong G_2(q)$ or $^3D_4(q)$ and $i$ be an involution in $G$, then $C_G(i)$ is isomorphic to $\sl_2(q)\circ_2 \sl_2(q)$ or $\sl_2(q)\circ_2 \sl_2(q^3)$, respectively. Note that the components of $C_G(i)$ are long and short root $\sl_2(q)$-subgroups in $G$. If $G \cong \/^3D_4(q)$, then the subgroup $K\cong \sl_2(q^3) \leqslant C_G(i)$ is a short root $\sl_2(q)$-subgroup in $G$.
If $K$ is any component of $C_G(i)$ and $z$ be the unique involution in $Z(K)$, then $i=z$ and $K = K^g$ for any $g \in C_G(z)$. Hence Algorithm \ref{long-sl2q} returns that $K$ is a long root $\sl_2(q)$-subgroup. Therefore, in these cases we use the following algorithm.

\begin{algo}\label{g2and3d4} \textsc{``Construction of a long root $\sl_2(q)$ in $G_2(q)$ and $^3D_4(q)$''}
\begin{itemize}
\item[] Input: A black-box group $G$ isomorphic to $G_2(q)$ or $^3D_4(q)$, and the size of the field $q$.
\item[] Output: A black-box group $K$ which is a long root $\sl_2(q)$-subgroup in $G$.
\end{itemize}
\end{algo}

\noindent\textsc{Description of the algorithm:}\\
We first distinguish $^3D_4(q)$ from $G_2(q)$. Note that $|G_2(q)| = q^6(q^6-1)(q^2-1)$ and $|\/^3D_4(q)| = q^{12}(q^8+q^4+1)(q^6-1)(q^2-1)$. Hence, if we find an element $g\in G$ such that $g^m\neq 1$, where $m=q^6(q^6-1)(q^2-1)$, then we deduce that $G\cong \/^3D_4(q)$. Recall that the proportion of such elements corresponds to the sizes of the conjugacy classes of tori, which depend only on the Weyl group, so it is bounded from below by a constant. 

Assume that $G\cong  G_2(q)$. Let $i\in G$ be an involution, then $C_G(i) = L_1 \circ_2 L_2$, where $L_k \cong \sl_2(q)$, $k=1,2$. Note that $L_1$ and $L_2$ are short and long root $\sl_2(q)$-subgroups in $G$. Assume that $L_1$ (resp. $L_2$) is a short (resp. long) root $\sl_2(q)$-subgroup in $G$. Now let $j$ be an involution in $C_G(i)$ which does not centralize $L_1$ and $L_2$. We have $C_G(j) \cong \sl_2(q) \circ_2 \sl_2(q)$ since there is only one conjugacy class of involutions in $G$. Let $K_1$ and $K_2$ be short and long root $\sl_2(q)$-subgroups in $C_G(j)$, respectively. Then it is easy to see that all the pairs $\{L_s, K_t \}$, $s,t = 1,2$, generate $G_2(q)$ except $\langle L_2, K_2 \rangle \cong \sl_3(q)$ or $\su_3(q)$ \cite[Theorem 2.1]{liebeck94.293}. Now we need to distinguish $G_2(q)$ from $\sl_3^\eps(q)$ to recognize the long root $\sl_2(q)$-subgroups. To do this, we compute the number of components in a centralizer of an involution by following the arguments in Remark \ref{comm-sl2}. It is clear that the number of components in a centralizer of an involution in $\sl_3^\eps(q)$ is one, whereas there are two components in a centralizer of an involution in $G_2(q)$.

Let $G\cong \/^3 D_4(q)$ and $i \in G$ be an involution. Then $C_G(i) \cong \sl_2(q) \circ_2 \sl_2(q^3)$ where $\sl_2(q)$ corresponds to a long root and $\sl_2(q^3)$ corresponds to a short root $\sl_2(q)$-subgroup in $G$. To construct $\sl_2(q)$, let $S$ be a set of generators for $C_G(i)$. Setting $m=q(q-1)(q+1)$, we consider 
$$S^m = \{ g^m \mid g \in S\}.$$ 
Now $L=\langle S^m \rangle \cong \sl_2(q^3)$ and we construct an element $g \in C_{C_G(i)}(L)$. Note that $(q-1,q^3+1) = (q+1,q^3-1) =2$. Therefore we look for elements $g \in C_G(i)$ satisfying one of the conditions $g^{(q-1)(q^3+1)}=1$ and $o(g^{q^3+1})>2$ or $g^{(q+1)(q^3-1)}=1$ and $o(g^{q^3-1})>2$. Let $g=(g_1,g_2)\in C_G(i)$, $g^{(q-1)(q^3+1)}=1$, where $g_1^{q-1}=1$, $g_2^{q^3+1}=1$ and $g^{q^3+1}$ is non-central in $C_G(i)$. By Lemma \ref{pseudo-prod-sl2}, we find such elements with probability at least $1/6$ since the probability of finding a non-central element $g_1 \in \sl_2(q)$ of order dividing $q-1$ is at least $1/2$ and the probability of finding an element $g_2 \in \sl_2(q^3)$ of order dividing $q^3+1$ is at least $1/3$. Now, setting $h=g^{q^3+1}$ we have $\langle h^{C_G(i)} \rangle ^\prime = \sl_2(q)$ which corresponds to a long root $\sl_2(q)$-subgroup in $G$.

\subsection{The Main Algorithm}
We can now present Algorithm \ref{main-simple-alg}.

\noindent\textsc{Description of Algorithm \ref{main-simple-alg}:}
\begin{itemize}
\item[1.] Run Algorithm \ref{prod-sl2q} to construct a commuting product $L$ of the subgroups $\ppsl_2(q^l)$ for various $l\geqslant 1$.
\item[2.] Run Algorithm \ref{sl2q} to construct a subgroup $K\cong \sl_2(q^l)$ for some $l\geqslant 1$, which is a component of $L$.
\begin{itemize}
\item [2.1] If the number of recursive steps in Algorithm \ref{prod-sl2q} is more than one, then we go to Step 3. Otherwise, we check if $G \cong G_2(q)$ or $^3D_4(q)$.
\item[2.2] If $G$ is isomorphic to $G_2(q)$ or $^3D_4(q)$, then we run Algorithm \ref{g2and3d4} to construct a long root $\sl_2(q)$-subgroup. Otherwise, we go to Step 3. 
\end{itemize}
\item[3.] Use Algorithm \ref{long-sl2q} to check if $K$ is a long root $\sl_2(q)$-subgroup. If not, then construct another component of $L$ and check if it is a  long root $\sl_2(q)$-subgroup. If none of  the components of $L$ is a long root $\sl_2(q)$-subgroup, then we start the procedure from the beginning without going through subroutines 2.1 and 2.2 in Step 2.
\end{itemize}

{\bf Step 1:} We run Algorithm \ref{prod-sl2q} as described in Subsection \ref{find-prosl2}. If Algorithm \ref{prod-sl2q} fails to succeed in constructing a commuting product of subgroups $\ppsl_2(q)$ for various $q$, for example, it may return the identity group, then the construction of the component(s) of the centralizers of involutions in one of the recursive steps has failed. In this case, we produce more generators for the centralizers of involutions which are already constructed in the recursive steps. Thus we assume that Algorithm \ref{prod-sl2q} returns a commuting product $L$ of subgroups $\ppsl_2(q^l)$ in $G$,  where $l\geqslant 1$ may vary.

{\bf Step 2:} If the number of recursive steps in Algorithm \ref{prod-sl2q} is more than one, then we conclude that $G$ is not isomorphic to $G_2(q)$ or $^3D_4(q)$, and we go to Step 3.

Now assume that the number of recursive steps in Algorithm \ref{prod-sl2q} is one. By Table \ref{cent-table}, this can only happen in the following groups:
\[
\psl_3^\eps(q), \psl_4^\eps(q),\psp_4(q), \psizo_7(q), \pso_8^\pm(q), G_2(q),\/^3\/D_4(q).
\]
Since we have special subroutines to construct a long root $\sl_2(q)$-subgroup in $G_2(q)$ and $^3D_4(q)$, we need to distinguish them from the other groups listed above. To do this, we first compute the size, $q$, of the field. Let $K\cong \sl_2(q^l) \leqslant G$, where $l\geqslant 1$, be a subgroup constructed by Algorithm \ref{sl2q}. Now we check if Algorithm \ref{long-sl2q} returns that $K$ is a long root $\sl_2(q)$-subgroup in $G$. If not, then $G\ncong G_2(q)$ or $^3D_4(q)$, and we go to Step 3. If Algorithm \ref{long-sl2q} returns that $K$ is a long root $\sl_2(q)$-subgroup, then we conclude that $l=1$ except possibly $G\cong \/^3D_4(q)$ and $K\cong \sl_2(q^3)$. 

Recall that we compute the size of the field in Algorithm \ref{long-sl2q} by using \cite[Algorithm 5.5]{suko01}. Assume that we compute the size of the field as $q_0=p^k$ for some $k\geqslant 1$. Since $p$ is given as an input, we can compute $k$. If $k$ is not divisible by $3$, then clearly $q=q_0$. Now assume that $k$ is divisible by $3$. Let $q_1=p^{k/3}$ and 
$m=|^3D_4(q_1)|=q_1^{12}(q_1^2-1)(q_1^6-1)(q_1^8+q_1^4+1)$. Then we look for an element $g \in G$ such that $g^m \neq 1$. If we can not find such an element, then we conclude  that $q=q_1$ and $G\cong\/ ^3D_4(q)$, and we run Algorithm \ref{g2and3d4} to construct long root $\sl_2(q)$-subgroup. Observe that if $G\ncong\/ ^3D_4(q)$, then $q=q_1^3$, and we can find an element $g$ in all the groups listed above such that $g^m \neq 1$ with high probability. Hence we have computed $q$.

Now the groups $\psl_4(q)$, $\psp_4(q)$, $\pso_8^\pm(q)$ and $\psizo_7(q)$ have elements of order dividing $q^4-1$ but not $q^6-1$ so $G_2(q)$ and $^3D_4(q)$ can be distinguished from these groups since $G_2(q)$ and $^3D_4(q)$ do not have such elements. Similarly, we can distinguish $^3D_4(q)$ from $G_2(q)$ and $\psl_3^\eps(q)$ since $^3D_4(q)$ has elements of order dividing $q^6(q^8+q^4+1)$ but not $q^6(q^6-1)$, and $G_2(q)$ and $\psl_3^\eps(q)$ do not have such elements. Finally, we use the same arguments in Algortihm \ref{g2and3d4} to distinguish $G_2(q)$ from $\psl_3^\eps(q)$.

Alternatively, assuming the existence of an order oracle, $G_2(q)$ and $^3D_4(q)$ can be distinguished from the other groups listed above by using statistics of element orders \cite{babai02.383}. 

If $G \cong G_2(q)$ or $^3D_4(q)$ then we run Algorithm \ref{g2and3d4} to construct a long root $\sl_2(q)$-subgroup in $G$. Otherwise we go to Step 3.

{\bf Step 3:} We run Algorithm \ref{long-sl2q} to check whether $K$ is a long root $\sl_2(q)$-subgroup in $G$. If $K$ is not a long root $\sl_2(q)$-subgroup, then we construct another component of $L$ by following the arguments in Remark \ref{comm-sl2}.

Notice that if $K$ is not a long root $\sl_2(q)$-subgroup, then we do not need to go through the subroutines 2.1 and 2.2.

\section{Estimates}\label{estimates}
In this section we estimate the probability of producing an involution $i={\rm{i}}(g)$ from a random element $g \in G$ where the recursive construction of centralizers of involutions applied to $C_G(i)$ returns a commuting product of long root $\sl_2(q)$-subgroups in $G$. The estimates in some cases are very crude and on the cautious side, the actual probabilities of success are much higher.

Let $n=2^km$, $m$ odd. Then the number $k$ is called the \textit{$2$-height} of $n$.

Recall that we need an element of even order to produce involutions in $G$ and the proportion of elements of even order is at least $1/4$ by \cite[Corollary~5.3]{isaacs95.139}.

Assume that $G\cong \psl_{n}(q)$ and $i$ is an involution of type $t_{n/2}^\prime$. 
Then $C_G(i)^{\prime \prime} \cong \frac{1}{(n/2,q-1)}\sl_{n/2}(q^2).$
If an involution of type $t_{n/2}^\prime$ is constructed in one of the recursive steps, then we obtain a subgroup $\sl_2(q^{2^k})$ for some $k>1$ which is not a long root $\sl_2(q)$-subgroup in $G$. Note that an involution of type $t_{n/2}^\prime$ belongs only in a maximal twisted torus $T$, twisted by the longest cycle $w$ in the Weyl group $W \cong {\rm{Sym}}(n)$. Therefore $i={\rm{i}}(g)$ where $g$ is an element from a maximal twisted torus. The number of maximal twisted tori in $G$ is $|G:N_G(T)|$. Moreover, $|N_G(T):T| = |C_W(w)|$ by Equation (\ref{toriat-q}) in Section \ref{maximal_tori} and $|C_W(w)|=n$ since $w$ is a cycle of length $n$ in $W$. Hence, there are at most 
\[
\frac{|G|}{|N_G(T)|}= \frac{|G|}{|C_W(w)|\cdot |T|}=\frac{|G|}{n|T|}
\]
elements which produce involutions of type $t_{n/2}^\prime$. Therefore the probability of producing an involution of type $t_{n/2}^\prime$ is at most $1/n$. The other types of involutions have centralizers of the form  $L\cong \sl_k(q) \circ \sl_l(q)$ in which case producing an involution $i={\rm{i}}(g)$ of type $t_{k/2}^\prime$ or $t_{l/2}^\prime$ in the components is close to $0$ as we must have an element $g \in L$ whose components belong to the same type of torus and have the same 2-height. Therefore if we have an involution of type different than $t_{n/2}^\prime$ in the first step, which is of probability at least $1/4(1-1/n)$, then Algorithm \ref{prod-sl2q} returns a commuting product of long root $\sl_2(q)$-subgroups with probability close to 1. The case where $G\cong \psu_n(q)$ is analogous. 

Assume that $G\cong  \psp_{2n}(q)$ and $q \equiv 1\, ({\rm mod} \, 4)$. Let $i$ be an involution of type $t_n$, then $C_G(i)^{\prime \prime} \cong \frac{1}{(2,n)}\sl_n(q)$. If we construct an involution of type $t_n$ in one of the recursive steps, then Algorithms \ref{prod-sl2q} and \ref{sl2q} do not return a long root $\sl_2(q)$-subgroup. To see this, let $G\cong \sp_{2n}(q)$, $V$ be the natural module for $G$ and $t$ be the preimage of the involution $i$. Then $C_G(i)$ is a maximal subgroup of $G$ which leaves invariant some totally isotropic subspace of $V$. Therefore if $K\cong \sl_2(q)$ is the subgroup obtained by recursive construction of centralizers of involutions applied to $C_G(i)$ and $W$ is the subspace on which $K$ acts, then the  symplectic form on $W$ is degenerate. Now, it is easy to see that the involutions of type $t_n$ belong to some maximal twisted torus $T\leqslant G$. The number of maximal twisted tori in $G$ is $|G:N_G(T)|$ and $|N_G(T):T| = |C_W(w)|=2n$ since $T$ corresponds to a longest cycle $w$ in the Weyl group $W=Z_2 \wr {\rm Sym}(n)$, see Lemma 2.3 in \cite{altseimer01.1}. Hence, there are at most 
$$\frac{|G|}{2n|T|}$$ 
elements which can produce involutions of type $t_n$. Therefore the probability of producing an involution of type $t_n$ is at most $1/2n$. Hence, we can produce an involution which is not of type $t_n$ with probability at least $1/4(1-1/2n)$. If $q \equiv -1\, ({\rm mod} \, 4)$, then $C_G(i)^{\prime \prime} = \frac{1}{(2,n)}\su_n(q)$, where $i$ is an involution of type $t_n$ and the same arguments apply to obtain the same estimate. As above, the other types of involutions have centralizers of the form $\sp_{2k}(q) \circ \sp_{2n-2k}(q)$ and so the probability of producing an involution of type $t_k$ or $t_{n-k}$ from random elements in these components is close to $0$ by the same arguments.

Assume that $G\cong  \pso_{2n}^\varepsilon(q)$, $n\geqslant 4$. If we construct $\frac{1}{2}\sl_n^\varepsilon(q)$ as a centralizer of an involution, then a lower bound for the probability of  constructing a long root $\sl_2(q)$-subgroup follows from the similar estimate for $\psl_n(q)$. The desired involution is of type $t_n$ or $t_n^\prime$. Again these involutions belong to a maximal twisted tori. By using the same ideas above the probability of obtaining such an involution is at least $1/2n$, and the overall probability is $1/4(1-1/n)(1/2n)$.

Assume that $G\cong \psizo_{2n+1}(q)$, $n \geqslant 3$. Let $i$ be an involution of type $t_n$. Then $C_G(i)^{\prime\prime} \cong \psizo_{2n}^\eps(q)$, where $q^n \equiv \eps \, ({\rm mod}\, 4)$.  Note that only maximal twisted tori whose orders are $1/2(q^n\pm 1)$ contain involutions of type $t_n$. Since $\psizo_{2n+1}(q)$ and $\psp_{2n}(q)$ have the same Weyl groups, the estimate for the construction of an involution of type $t_n$ in the first recursive step of Algorithm \ref{prod-sl2q} is the same as in the case of $\psp_{2n}(q)$ which is $1/2n$. Therefore a crude estimate in this case follows from an estimate in $\psizo_{2n}^\eps(q)$. 

Assume that $G\cong F_4(q)$. If we have an involution $i$ of type $t_4$, then $C_G(i)^{\prime \prime} \cong {\rm Spin}_9(q)$ and the estimate for constructing a long root $\sl_2(q)$-subgroup follows from the estimate for $\psizo_9(q)$. An involution of type $t_4$ belongs to a torus $T$, where $T$ corresponds to an element $w \in W$ with $|C_W(w)|=8$ \cite[Table 4]{carter70.g1}.

Let $G\cong  E_8(q)$ or $E_7(q)$. Then the semisimple socles of the centralizers of involutions are either central products of classical groups or contain $E_7(q)$ or $E_6(q)$, respectively. In the case of central products of classical groups, we refer to the above estimates. Similarly, we reduce the estimates for the groups $E_7(q)$ to the estimates in $E_6(q)$. The estimates for the groups $E_6(q)$ and $^2 E_6(q)$ can be computed again from the estimates for classical groups as the semisimple socles of the centralizers of involutions are central products of classical groups, see Table \ref{cent-table}.

Recall that the construction of a long root $\sl_2(q)$-subgroup in $G_2(q)$ and $^3D_4(q)$ is not a recursive procedure and it can be constructed without any difficulty by the routine presented in Section \ref{sec-g2and3d4}.

\section{Recognition of the $p$-core}\label{p-core}
The aim of this section is to prove Theorem \ref{maincor1}.  The algorithm is primarily focused on the construction of the centralizers of classical involutions. We use the function $\zeta_1$ to generate the centralizers of classical involutions. Therefore it is enough to obtain a lower bound for the probability that $i i^g$ has odd order only for classical involutions $i \in G$.

\begin{theorem}\label{zeta1-lemma}
Let $G$ be a finite simple classical group over a field of odd characteristic $p$ and $i \in G$ be a classical involution. Then the product $i i^g$ has odd order with probability bounded from below by a constant which does not depend on $G$.
\end{theorem}
\proof Let $i \in G$ be a classical involution. Then it belongs to some long root $\sl_2(q)$-subgroup $K\leqslant G$. Therefore, the product $ii^g$ belongs to the subgroup $L=\langle K, K^g \rangle$. We have shown in Section \ref{pairs} that the subgroup $L$ has a given structure depending on $G$ with probability at least $1-O(1/q)$. Hence, it is enough to find the probability that a product of two conjugate classical involutions in $L$ has odd order. %The arguments for computing the estimates are same as in the proof of Lemma 2.9 in \cite{altseimer01.1}. 

Consider the map 
\begin{eqnarray*}
\varphi: i^L \times i^L & \rightarrow & L \\
(i^g, i^h) & \mapsto & i^gi^h.
\end{eqnarray*}
Take a torus $T\leqslant L$ inverted by $j=i^{g^\prime}$ for some $g^\prime \in L$, that is, satisfying $t^j =t^{-1}$ for all $t \in T$. Let $x \in T$ be an element of odd order. Then there exists $h \in \langle x \rangle$ such that $h^2 =x$. Now 
$$jj^h=jh^{-1}jh=hh=x$$ 
since $j$ inverts $T$. Hence elements of odd order in $T$ are in the image of $\varphi$. Let $x\in T$ be a regular element, that is, $C_L(x) =T$, which has odd order. Then $|\varphi^{-1}(x)| = |T|$ since  
$j^tj^{ht} = (jj^h)^t=(hh)^t=x^t=x$ for any $t \in T$. Let $S$ be the set of regular elements in $T$ which are of odd order. Let $R$ be the set of all regular elements in $L$ whose elements are conjugate to elements in $S$, then
$$|R| =  |L: N_L(T)||S|.$$
Now observe that 
$$|\varphi^{-1}(R)| \, \geqslant \, |R||S|, $$
and $$|i^L \times i^L| = \frac{|L|^2}{|C_L(i)|^2}.$$
Therefore the proportion of pairs of involutions which are mapped to $R$ is 
$$ \frac{| \varphi^{-1}(R) |}{| i^L \times i^L |}  \geqslant \frac{| R | | S| | C_L(i)|^2}{| L |^2} = \frac{ | S|^2 | C_L(i)|^2}{| N_L(T)| | L|}.$$

Now, we shall find a lower bound for this quotient in all finite simple classical groups.

Let $G\cong \psl_n(q)$, $n\geqslant 5$. Then by Lemma \ref{psl}, $L\cong \sl_4(q)$ with probability at least $1-1/q^{n-3}>1/2$. Observe that a classical involution $i \in L$ inverts a cyclic torus $H\leqslant L$ of order $q^2+1$ and $(q^2+1)/2$ is odd. Observe also that $H$ is uniquely contained in a maximal cyclic torus $T \leqslant L$ of order $(q+1)(q^2+1)= (q^4-1)/(q-1)$. Note that the maximal torus $T$ in $L$ corresponds to a $4$-cycle in the Weyl group of $L$ which is $Sym(4)$ in this case. Hence, $|R| = \frac{|L|}{8 (q+1)}$, $|S| = (q^2+1)/2$, $|C_L(i)| = q^2(q^2-1)^2(q-1)$ and $|L| = q^6(q^2-1)(q^3-1)(q^4-1)$. After a simple rearrangement 
$$\frac{| \varphi^{-1}(R) |}{| i^L \times i^L |} \geqslant \frac{(q-1)^3(q+1)}{16q^2(q^2+q+1)} \geqslant \frac{1}{32}.$$ 
Hence, $ii^g \in \psl_n(q)$ is of odd order with probability at least $1/64$. If $G\cong \ppsl_n(q)$ with $n\leqslant 4$, then by Lemma \ref{others}, $L=G$ with probability at least $1-1/q$, and by the same computations we obtain a similar result.
 
Let $G\cong \psp_{2n}(q)$, $n\geqslant 3$. Then $L \cong \sp_4(q)$ with probability at least $1-1/q^{2n-2}>1/2$ by Lemma \ref{psp}. The classical involutions invert a tori of order $q \pm 1$. Let $T\leqslant L$ be a torus of order $q-1$ inverted by a classical involution $i \in L$ and $q \equiv -1\, ({\rm mod}\, 4)$. Then $(q-1)/2$ is odd. Observe that $C_L(T) \cong \gl_2(q)$, $| N_L(T)/C_L(T) | = 2$ and $C_L(i) \cong \sl_2(q) \times \sl_2(q)$. Now after a simple rearrangement, we have
$$\frac{| \varphi^{-1}(R) |} {| i^L \times i^L |} \geqslant \frac{1}{4}\frac{(q-1)^2(q+1)}{q(q^2+1)}\geqslant 1/16.$$
Hence, $ii^g \in \psp_{2n}(q)$ is of odd order with probability at least $1/32$. If $q \equiv 1\, ({\rm mod} \,4)$, then we consider tori of order $q+1$. If $G\cong \psp_4(q)$, then by Lemma \ref{others}, $L=G$ with probability at least $1-1/q$ and by the same arguments we obtain a similar result.

Let $G\cong \psizo_8^{-}(q)$, then $L=G$ with probability at least $1-1/q>1/2$ by Lemma \ref{others}. Observe that a classical involution $i \in G$ inverts a maximal torus $T$ of order $(q^4+1)/2$ which is odd. Now, $|C_G(i)| = \frac{1}{4}q^4(q-1)^3(q+1)^3(q^2+1)$, $|S| = (q^4+1)/2$, $|N_G(T)/T| = 12$. Hence,
$$\frac{| \varphi^{-1}(R) |} {| i^L \times i^L |} \geqslant \frac{1}{384}\frac{(q-1)^4(q+1)^4(q^2+1)}{q^4(q^6-1)} \geqslant \frac{1}{4\cdot 384}=\frac{1}{1536}.$$
Hence, $ii^g \in \psizo_8^-(q)$ is of odd order with probability at least $1/(2\cdot 1536)$. 

Assume now that $G\cong  \psizo_8^+(q)$ or $\psizo_n^\eps(q)$, $n\geqslant 9$. Then $L \cong \psizo_8^+(q)$ with probability at least $1-1/q>1/2$ by Lemma \ref{psizo} and \ref{others}. Observe that $L$ contains a subgroup of the form $N\cong \psizo_4^-(q) \times \psizo_4^-(q)$, and there is an involution $i \in \psizo_4^-(q)$ which inverts a torus of order $(q^2+1)/2$ in $\psizo_4^-(q)$. Notice that $i$ is necessarily an involution of type $t_1$ in $L$. Now the involution $j=(i,i) \in N$ inverts a torus $T$ of order $(q^2+1)^2/4$ which is odd. Since $j$ is a product of two commuting involutions of type $t_1$, it is of type $t_2$ in $L$ and therefore it is a classical involution. Hence,
$$| S|=(q^2+1)^2/4,$$
$$| C_L(i) | = 4 | \psizo_4^+(q)|^2 = q^4(q-1)^4(q+1)^4,  $$
$$| L | = q^{12}(q^4-1)(q^2-1)(q^4-1)(q^6-1), $$
$$| N_L(T)| = 8(q^2+1)^2,$$
and
$$ \frac{| \varphi^{-1}(R) |}{| i^L \times i^L |}  \geqslant 
   \frac{q^8(q-1)^8(q+1)^8(q^2+1)^4}{128(q^2+1)^2(q^{12}(q^4-1)(q^2-1)(q^4-1)(q^6-1))}.$$
After a simple rearrangement we have
$$ \frac{| \varphi^{-1}(R) |}{| i^L \times i^L |}  \geqslant 
   \frac{(q^2-1)^5}{128q^4(q^6-1)} \geqslant \frac{1}{128\cdot 6} = \frac{1}{768}. $$
\hfill $\Box$

The following lemma is crucial for our algorithm.

\begin{lemma}\label{invo-cent}
Let $X$ be a finite group, $i \in X$ an involution and $Q=O_p(X)$. If $C_Q(i)=1$ then $[i,x] \in C_X(Q)$ for all $x \in X$. 
\end{lemma}
\proof Notice that if $C_Q(i)=1$, then $i$ inverts $Q$, that is, $x^i=x^{-1}$ for any $x \in Q$ and $Q$ is abelian.  It is now easy to see that $[i,x] \in QC_X(Q)=C_X(Q)$ for all $x \in X$.  \hfill $\Box$

First, we present our algorithm in the base case, that is, $X/O_p(X)  \cong\psl_2(q)$. Let $i \in X$ be an involution and $Q = O_p(X)\neq 1$. We can assume that $C_Q(i) \neq 1$ by Lemma \ref{invo-cent} since otherwise random elements in $X$ power up to $p$-elements in $Q$ with high probability. Now $O_p(C_X(i)) \neq 1$ and $C_X(i)/O_p(C_X(i))$ is isomorphic to a dihedral group of order $q\pm 1$. Let $Q_1 = O_p(C_X(i))$. If $O_p(C_X(i)^\prime) = 1$, then random elements in $C_X(i)$ have orders which are multiples of $p$ and we can find a $p$-element in $Q_1$ by raising a random element in $C_X(i)$ to the power $q \pm 1$. Hence, we can assume that $O_p(C_X(i)^\prime) \neq 1$. Now $C_X(i)^\prime/O_p(C_X(i)^\prime)$ is isomorphic to a cyclic group of order $(q \pm 1)/2$ in which case we take the power $(q\pm 1)/2$ of random elements in $C_X(i)^\prime$ to produce $p$-elements in $O_p(C_X(i)^\prime)$ and we are done. Our approach in the general case is to reduce the problem to this base case in all finite simple classical groups. 

\begin{algo}{\rm (cf. \cite[Algorithm 5.8]{suko01})} \textsc{``Recognition of the $p$-core''}\label{pcore-alg}
\begin{itemize}
\item[] Input: A black-box group $X$ with the property that $X/O_p(X)$ is a finite simple classical group of odd characteristic $p$.
\item[] Output: Either the statement ``$O_p(X) \neq 1$'' together with a non-trivial $p$-element $g\in O_p(X)$ or the statement ``Possibly, the $p$-core is trivial''. 
The answer $``O_p(X)\neq 1''$ is always correct whereas the negative answer ``Possibly, the $p$-core is trivial'' may be false with some probability of error.
\end{itemize}
\end{algo}

\noindent\textsc{Description of the algorithm:}
\begin{itemize}
\item[1.] Check whether random search works in $X$. 
%\item[2.] Construct an involution $i \in X$ and check whether the elements $[i,x]$ power up to $p$-elements in $O_p(X)$. If we can not find a $p$-element, then we go to the next step.
%\item[3.] Construct $C_X(i)$ by using $\zeta_0^i\sqcup\zeta_1^i$, and the subgroups $C_X(i)^\prime$ and $C_X(i)^{\prime \prime}$. Check whether random search works in these subgroups.
%\item[4.] Apply Algorithm \ref{main-simple-alg} to the subgroup $C_X(i)^{\prime \prime}$ to construct a subgroup $K \leqslant X $ where $K/O_p(K)$ is a long root $\sl_2(q)$-subgroup in $X/O_p(X)$ and check whether $O_p(K) \neq 1$.
%\item[5.] Construct $C_X(z)$ by using $\zeta_1^z$, and the subgroups $C_X(z)^{\prime}$ and $C_X(z)^{\prime \prime}$ where $z\in Z(K)$. Check whether random search works in these subgroups.
%\item[6.] Construct the subgroups $K_1, L_1$ where $C_X(z)^{\prime \prime}=K_1L_1$ and $[K_1,L_1]=1$. Here, $K_1/O_p(K_1)$ is a long root $\sl_2(q)$-subgroup in $X/O_p(X)$. Check whether random search works in $K_1$ and $L_1$.
%\item[7.] Set $X=L_1$ and go to Step 1.
\item[2.] Construct a subgroup $K\leqslant X$ where $K/O_p(K)$ is a long root $\sl_2(q)$-subgroup in $X/O_p(X)$ and check whether $O_p(K)\neq 1$.
\item[3.] Construct $C_X(i)$ by using $\zeta_1^i$, and the subgroups $C_X(i)^\prime$ and $C_X(i)^{\prime \prime}$ where $i \in Z(K)$. Check whether random search works in these subgroups.
\item[4.] Construct the subgroups $K_1, L_1$ where $C_X(i)^{\prime \prime}=K_1L_1$ and $[K_1,L_1]=1$. Here, $K_1/O_p(K_1)$ is a long root $\sl_2(q)$-subgroup in $X/O_p(X)$. Check if $O_p(K_1) \neq 1$. If $O_p(K_1)=1$, then go to the next step.
\item[5.] Set $X=L_1$ and go to Step 1.
\end{itemize}

{\bf Step 1:}  First, we take random elements and check whether they power upto $p$-elements. To do this, we compute the natural number $m$, where $E=mp^k$, $(m,p)=1$, is the exponent of  $X$ given as an input. We can check whether a $p$-element $g \in X$ belongs to the $p$-core $O_p(X)$ in the following way: We construct the normal closure $P=\langle g \rangle^ X$ and then check the solvability of this subgroup. By \cite[Theorem 1.5]{babai95.296}, we can construct the normal closure of a subgroup and decide the solvability of a given black-box group. Now, if $P$ is not solvable, then $g \notin O_p(X)$.

In practice, one can check whether a $p$-element $g\in X$ belongs to $O_p(X)$ in the following way. Produce a random element $h \in X$ and construct the subgroup $P = \langle g, g^h \rangle$. If $g \notin O_p(X)$, then $P$ is not a $p$-group with high probability.

It is proved in \cite{babai01.39} that if $X/O_p(X)$ is a finite simple unisingular group of Lie type, then random elements power up to a $p$-element in $O_p(X)$ with high probability.

The classification of the finite simple unisingular groups is as follows.

\begin{theorem}{\rm (\cite[Theorem 1.3]{guralnick03.271})}
A finite simple group $G$ of Lie type of characteristic $p$ is unisingular if and only if $G$ is one of the following:
\begin{itemize}
\item[(i)] $\psl_n^\eps(p)$ with $n \mid (p-\eps)$;
\item[(ii)] $\psizo_{2n+1}(p), \, \psp_{2n}(p)$ with $p$ odd;
\item[(iii)] $\pso_{2n}^\eps(p)$ with $p$ odd, $\eps=(-1)^{n(p-1)/2}$;
\item[(iv)] $^2G_2(q)$, $F_4(q)$, $^2F_4(q)$, $E_8(q)$ with $q$ arbitrary;
\item[(v)] $G_2(q)$ with $q$ odd;
\item[(vi)] $E_6^\eps(p)$ with $3 \mid (p-\eps)$;
\item[(vii)] $E_7(p)$ with $p$ odd.
\end{itemize}
\end{theorem}

{\bf Step 2:} We use Algorithm \ref{main-simple-alg} to construct a subgroup $K$ where $K/O_p(K)$ is a long root $\sl_2(q)$-subgroup in $X/O_p(X)$. In each recursive step in the construction of $K$, we check whether random search works as in Step 1. Note that if $C_{O_p(X)}(i) =1$, where $i\in X$ is an involution, then $i$ inverts $O_p(X)$ by Lemma \ref{invo-cent}. Moreover, the elements of the form $[i,x]$ power up to a $p$-element in $O_p(X)$ with high probability. Therefore, before constructing centralizer of an involution in each recursive step, we also check whether that involution inverts the $p$-core of the subgroup constructed in the previous recursion. 

%we also check whether the elements $[i,x]$ power up to $p$-elements in $O_p(X)$ for random elements $x \in X$.

%whether the involution inverts the $p$-core. That is, we check whether the elements $[i,x]$ power up to $p$-elements in $O_p(X)$, where $i\in X$ is an involution and $x\in X$ is a random element. 

Notice that we may have $O_p(K)=1$ even though $O_p(X) \neq 1$. We check whether $O_p(K) \neq 1$ as above. If we can not find a $p$-element then we go to the next step.

{\bf Step 3:} We use $\zeta_1^i$ to construct $C_X(i)$. The derived subgroups $C_X(i)^\prime$ and $C_X(i)^{\prime \prime}$ can be constructed by an algorithm in \cite{babai95.296}. As discussed in Subsection \ref{concgi}, the map $\zeta_1^i$ can be used efficiently to generate $C_X(i)$ by Theorems \ref{dist} and \ref{zeta1-lemma}.

If $Q\neq 1$ and $C_Q(i) = 1$, then by Lemma \ref{invo-cent}, $[i,x] \in C_X(Q)$ for all $x \in X$ and these elements power up to $p$-elements in $Q$ with high probability. Therefore we assume that $C_Q(i) \neq 1$ which implies that $O_p(C_X(i)) \neq 1$. If $O_p(C_X(i)^\prime) = 1$ or $O_p(C_X(i)^{\prime\prime}) = 1$, then again random elements in $C_X(i)$ or $C_X(i)^\prime$ power up to $p$-elements in $Q$, respectively. Therefore we assume now that $O_p(C_X(i)^{\prime \prime}) \neq 1$ and go to the next step.

{\bf Step 4:} We construct the 2-components $K_1$ and $L_1$ of $C=C_X(i)^{\prime \prime}$; compare with \cite[Section 11]{leedham-green09.833} for similar computations. By Theorem  \ref{long-root-sl2} and Corollary \ref{semisocle2}, $C=K_1L_1$, where $[K_1,L_1]=1$. Here, $K_1/O_p(K_1)$ is a long root $\sl_2(q)$-subgroup in $X/O_p(X)$ and the structure of the subgroup $L_1$ can be read from Table \ref{long-root}. Notice that $K\leqslant K_1$.

We consider a maximal twisted torus $T\leqslant L_1$. Recall that the maximal twisted tori are conjugate in $L_1$ and the probability that a random element belongs to a torus conjugate to $T$ is $O(1/n)$ by Equation (\ref{toriat-q}). 

%\begin{center}
\begin{table}[h]
\caption{The orders of some maximal twisted tori in $L_1/O_p(L_1)$ \cite{carter70.g1}, see also \cite[\S 3.1]{suko01}.}\label{ord-max-twist}
\begin{tabular}{cc|cc|cc}
$X/O_p(X)$ &&$L_1/O_p(L_1)$                  &&&  $| T |$                         \\[.9ex] \hline
$\psl_n(q)$&&$\sl_{n-2}(q)$         &&&   $(q^{n-2}- 1)/(q-1)$ or $q^{n-3}-1$               \\[.8ex]\hline
$\psu_n(q)$&&$\su_{n-2}(q)$         &&&  $(q^{n-2}- (-1)^{n-2})/(q+1)$         \\[.9ex]\hline
$\psp_{2n}(q)$&&$\sp_{2n-2}(q)$        &&&  $(q^{n-1} + 1)/2$                     \\[.9ex]\hline
$\psizo_{2n+1}(q)$&&$\psizo_{2n-3}(q)$     &&&  $(q^{n-2} + 1)/2$                     \\[.9ex]\hline
$\pso_{2n}^+(q)$&&$\psizo_{2n-4}^+(q)$   &&&  $(q^{n-3}+1)(q+1)/2$                   \\[.9ex]\hline
$\pso_{2n}^-(q)$&&$\psizo_{2n-4}^-(q)$   &&&  $(q^{n-2}+1)/2$                        
\end{tabular}
\end{table}
%\end{center}

Let $X/O_p(X)\cong \psl_n(q)$ and $L_1/O_p(L_1)\cong \sl_{n-2}(q)$. Observe that $((q^{n-3}-1),q+1)=2$ if $n$ is even, and $((q^{n-2}-1)/(q-1),q+1)=1$ if $n$ is odd. Therefore we consider tori of order $(q^{n-3}-1)$ or $(q^{n-2}-1)/(q-1)$ if $n$ is even or odd, respectively. Provided that $n\geqslant 5$, the probability that an element in $L_1$ having an order dividing $q^{n-3}-1$ or $(q^{n-2}-1)/(q-1)$ is at least $1/(n-3)$ or $1/(n-2)$, respectively. Therefore with probability at least $1/(n-2)$, $h=g^{E/(q+1)^a}\in K_1$, where $E$ is an exponent for $X$ and $a$ is the biggest power of $(q+1)$ in $E$.  If $h$ is a central element in $C$ then we repeat this proces until we find a non-central element. If $h\in C$ is a non-central element, then it is clear that $\langle h^C \rangle ^\prime =K_1$.

In the rest of the classical groups, except when $L_1/O_p(L_1)\cong \su_{n-2}(q)$ and $n$ is even, we take a torus $T$ as in the third column of Table \ref{ord-max-twist}. Observe that we have $(| T |, q-1) = 2$. If $L_1/O_p(L_1)\cong \su_{n-2}(q)$ and $n$ is even, then we take a torus $T$ of order $(q^{n-3} +1)$ as $(q^{n-3}+1, q-1)=2$. Therefore, after $O(n)$ iterations, we can find an element $g \in C$ such that $h = g^{E/(q-1)^b}$ is a non-central element in $K_1$, where $b$ is the maximal power of $q-1$ in $E$. Thus $\langle h^C \rangle^\prime = K_1$. 

Now we check whether $K_1$ has non-trivial $p$-core as before. If $O_p(K_1) = 1$, then we construct $L_1$ by raising the power $q(q^2-1)$ of the elements in the generating set for $C$. If $L_1$ is a commuting product of subgroups $\ppsl_2(q)$, then we use the procedure in Remark \ref{comm-sl2}  to construct each $\ppsl_2(q)$ in $L$.

{\bf Step 5:} If $O_p(X) \neq 1$ and $O_p(K_1) = 1$, then $O_p(L_1) \neq 1$. Now, we set $X = L_1$ and go to Step 1. In this way we construct a list of subgroups $K_s \leqslant X$, $s \geqslant 1$, where $K_s/O_p(K_s)$ is centrally isomorphic to $\psl_2(q)$ for each $s$. If we fail to construct $p$-elements in all these subgroups, we conclude that $O_p(X) =1$.

\section{Implementation}\label{gap}
The algorithms in this paper were tested in GAP4 Version 4.4.10 \cite{gap}. The GAP code was not written for practical purposes, but checking for the theoretical justifications. The author decided not to optimize the programming but to check whether each branch of the algorithms worked in practice. Therefore, there is a lot room for an improvement of the code for practical purposes. The source of the code can be obtained from the author upon request.

Although we do not calculate the order of an element in our justifications of the algorithms, we use, for simplicity, the  \verb Order  function to calculate the order of the elements. We use the product replacement algorithm to construct random elements in a given group. Instead of using GAP functions \verb DerivedSubgroup  and \verb NormalClosure  in the code, we produce random commutators and random conjugates of the generators to generate the derived subgroup and the normal closure of a subgroup. %In our experiments we produce 50 generators.

Although the code is not written for practical purposes, it is worthwhile to state the performance of  Algorithm \ref{main-simple-alg} and  we present some experiments which were carried out on a 2.4 GHz Intel Core 2 Duo processor in Table \ref{main-imp}. The running time is the CPU time in seconds which is avaraged over ten runs.

\begin{table}[h]
\caption{The implementation of Algorithm \ref{main-simple-alg}}\label{main-imp}
\begin{tabular}{cccc|cccc|ccccc}
$G$ && Time &&&$G$ &&Time &&&$G$ &&Time \\ \hline
$\sl_6(3^5)$ &&       1.55 &&&  $\psizo_8^+(5^3)$    && 1.76    &&& $E_6(7)$    && 7.59\\
$\sl_6(3^{10})$ &&   4.86&&& $\psizo_8^+(5^6)$   &&  6.2  &&& $E_6(7^5)$  &&137.25\\
$\sl_{20}(3^5)$ &&   13.1&&&   $\psizo_{16}^+(5^3)$   &&  6.57  &&&$E_7(7)$    &&42.74\\
$\sl_{20}(3^{10})$&& 78.82&&&  $\psizo_{16}^+(5^6)$  &&  34.25  &&& $E_7(7^3)$   && 742.67  \\
$\su_6(3^2)$ && 1.54 &&& $\psizo_7(5^3)$ &&   1.86   &&& $E_8(7)$    && 1979.68    \\
$\su_6(3^5)$ && 3.53   &&& $\psizo_7(5^6)$ &&   5.19   &&&$F_4(7)$  && 6.74    \\
$\su_{16}(3^2)$ &&  10.35  &&& $\psizo_{15}(5^3)$ &&  5.58    &&&$F_4(7^5)$  && 100.78    \\
$\su_{16}(3^5)$ &&  34.43  &&& $\psizo_{15}(5^6)$ && 28.3    &&& $^2E_6(7)$  && 12.32     \\
$\sp_{6}(3^5)$  &&  1.52     &&&  $\psizo_{10}^-(5^3)$    &&  2.37       &&& $G_2(7)$         && 1.89       \\
  $\sp_6(3^{10})$   &&   4.54    &&&  $\psizo_{10}^-(5^6)$    && 9.85  &&&   $G_2(7^5)$       &&  14.81       \\
$\sp_{20}(3^5)$     && 10.37      &&&  $\psizo_{20}^-(5^3)$              &&  10.25      &&& $^3D_4(7)$         &&   5.98      \\
 %$\sp_{20}(3^{10})$   &&  63.4      &&& $\psizo_{20}^-(5^6)$              &&    102.3    &&&          &&         \\
\end{tabular}
\end{table}

The performance of Algorithm \ref{pcore-alg} depends on the structure of the $p$-core $O_p(X)$, that is, even if $X/O_p(X)$ is not a unisingular group, Algorithm \ref{pcore-alg} may return a $p$-element before it constructs a long root $\sl_2(q)$-subgroup. For example, consider the affine group $G=H\ltimes V$, where $H\cong \sl_n(q)$ and $V$ is the natural module for $H$. Then the central involution in $ H$ inverts every element of $V$ so, in this particular case, Algorithm \ref{pcore-alg} returns a $p$-element before it starts constructing centralizers of involutions. 

Assume now that $X = G\ltimes H$, $G\cong \sl_n(q)$, $q=p^k$, $k\geqslant 2$ and $H$ is a $p$-group where $X$ is given in the following way: 
$$X= \left[
         \begin{array}{cc}
          G    &  H   \\
           0   & \hat{G}   \\
          \end{array}
      \right]. $$ 
Here $G$ is embedded diagonally in $\sl_{2n}(q)$ and the matrix entries of $\hat{G}$ are increased by $p$. It is clear that $H= O_p(X)$ and the involutions in $X$ do not invert $H$, that is, $C_H(i) \neq 1$ for any involution $i\in X$. Assuming that $q$ is very large, a $p$-element can not be constructed in one of the recursive steps corresponding to Step 2 of Algorithm \ref{pcore-alg}. However, in this case, Step 2 of Algorithm \ref{pcore-alg} returns a subgroup $K$ where $K/O_p(K)$ corresponds to a long root $\sl_2(q)$-subgroup in $X/O_p(X)$ and $O_p(K) \neq 1$. Hence, we construct a $p$-element in $K$ as described in Section \ref{p-core}. In this case the performance of the algorithm is similar to the performance presented in Table \ref{main-imp}.

\section*{Acknowledgements}
I am very grateful to Alexandre Borovik for many stimulating conversations and invaluable comments during the preparation of this work which would not have appeared without him. I would like to thank to Ayse Berkman for her useful suggestions. I would also like to thank to Tuna Alt\i nel and Frank Wagner for the invitation to visit Universit\'{e} Claude Bernard Lyon-1 where part of this work was carried out. This research is partially supported by the Turkish Scientific Council (T\"{U}B\.{I}TAK), MATHLOGAPS, project no: 504029, when the author visited Universit\'{e} Claude Bernard Lyon-1 and ARC Federation Fellowship FF0776186. I would like to acknowledge both the University of Manchester and Universit\'{e} Claude Bernard Lyon-1 for providing an excellent study environment during my stays.

\end{document}